\documentclass[12pt]{article}
\usepackage{kz}
\usepackage{enumerate}
\usepackage[font=scriptsize,labelfont=bf]{caption}
\RequirePackage[colorlinks,linkcolor=blue,citecolor=blue,urlcolor=blue]{hyperref}

\newcommand{\blind}{1}

\begin{document}

\def\spacingset#1{\renewcommand{\baselinestretch}%
{#1}\small\normalsize} \spacingset{1}


\if1\blind
{
  \title{\textbf{BET on Independence}\thanks{This paper is dedicated to Professor Lawrence D. Brown with the author's deepest appreciation for his teaching and mentoring.}}
  \author{Kai Zhang\thanks{
    Kai Zhang is Assistant Professor (E-mail: zhangk@email.unc.edu), Department of Statistics and Operations Research, University of North Carolina at Chapel Hill, Chapel Hill, NC 27599. The author gratefully acknowledges the support through \textit{NSF DMS-1613112 and NSF IIS-1633212}}\hspace{.2cm}\\
    University of North Carolina at Chapel Hill}
  \maketitle
} \fi

\if0\blind
{
  \bigskip
  \bigskip
  \bigskip
  \begin{center}
    {\LARGE\bf BET on Independence}
\end{center}
  \medskip
} \fi

\bigskip
\begin{abstract}
We study the problem of nonparametric dependence detection. Many existing methods may suffer severe power loss due to non-uniform consistency, which we illustrate with a paradox. To avoid such power loss, we approach the nonparametric test of independence through the new framework of binary expansion statistics (BEStat) and binary expansion testing (BET), which examine dependence through a novel binary expansion filtration approximation of the copula. Through a Hadamard transform, we find that the symmetry statistics in the filtration are complete sufficient statistics for dependence. These statistics are also uncorrelated under the null. By utilizing symmetry statistics, the BET avoids the problem of non-uniform consistency and improves upon a wide class of commonly used methods (a) by achieving the minimax rate in sample size requirement for reliable power and (b) by providing clear interpretations of global relationships upon rejection of independence. The binary expansion approach also connects the symmetry statistics with the current computing system to facilitate efficient bitwise implementation. We illustrate the BET with a study of the distribution of stars in the night sky and with an exploratory data analysis of the TCGA breast cancer data.
\end{abstract}

\noindent%
\vfill

\newpage
\spacingset{1.68} 
\section{Introduction}\label{sec: intro}
Independence is one of the most foundational concepts in statistics. It is also one of the most common assumptions in statistical literature. Thus verifying independence is one of the most important testing problems. If we are not able to check this crucial condition, then we are ``betting on independence'' at the risk of losing the validity of our conclusions. In this paper, we study the dependence detection problem in a distribution-free setting, in which we do not make any assumption on the joint distribution. We focus on the test of independence between two continuous variables, though the approach can be generalized for more variables. Without loss of generality, we consider $n$ i.i.d. observations from the copula $(U,V)$ whose marginal distributions are uniform over $[0,1].$ This copula can be obtained by transformations with marginal cumulative distribution functions (CDF) when they are known. In this case, $U$ and $V$ are independent if and only if their joint distribution $\Pb_{(U,V)}$ is the bivariate uniform distribution over $[0,1]^2$, denoted by $\Pb_0$. We also study the case when the marginal CDFs are unknown. In this case, we can use the empirical CDFs, and the test is about the independence of observed ranks. The theory and procedures are shown to be similar.

Tests of independence have been extensively studied in statistics and information theory. One of the most classical parametric methods is based on the Pearson correlation, which can be interpreted as a measure of linear relationship. Classical results in  \cite{Renyi1959} connect correlation and independence. Recent tests based on robust versions of correlation include \cite{han2017biometrika}. Existing nonparametric testing procedures can be roughly categorized into three main classes:

(a) The CDF approach, which compares the joint CDF and the product of marginal CDFs: This pioneer approach includes variants of the Kolmogorov-Smirnov test such as \cite{hoeffding1948non} and \cite{romano1989}.

(b) The distance and kernel based approach, which can be regarded as a generalization of the correlation: One important recent development on dependence measures is the distance correlation \citep{szekely2007, szekely2009brownian}, which possesses the crucial property that a zero distance correlation implies independence. Tests based on sample versions of the distance correlation \citep{szekely2013distance,szekely2013energy} have since been popular methods. Other important methods include the generalized measures of correlation (GMC) by \cite{zheng2012generalized} and the Hilbert Schmidt independence criterion (HSIC) by \cite{gretton2007kernel, sejdinovic2013, pfister2016kernel} who study dependence through distances between embedding of distributions to reproducing kernel Hilbert spaces (RKHS).

(c) The binning approach, which generalizes the comparison of the joint density and the product of marginal ones: By discretizing $X$ and $Y$ into finite many categories, classical statistical or information theoretical methods such as the $\chi^2$ tests and Fisher's exact tests can be applied to study the dependence. \cite{miller1982maximally} studied the maximal $\chi^2$ statistic from forming $2 \times 2$ tables through partitions of data. \cite{reshef2011detecting, reshef2015empirical, reshef2015equitability} introduced the maximal information coefficient (MIC) by aggregating information from optimal partitions of the scatterplot for different partition sizes. This approach was further studied by the $k$-nearest neighbor mutual information (KNN-MI) approach as described in \cite{kraskov2004estimating,kinney2014equitability}. \cite{heller2012consistent, heller2016consistent, heller2016multivariate} studied optimal permutation tests over partitions to improve the power. \cite{filippi2015bayesian} took a Bayesian nonparametric approach to the partitions. \cite{wang2016generalized} considered a generalized $R^2$ to detect piecewise linear relationships, a compromise between the distance approach and the binning approach that takes advantages of both. A very recent paper on Fisher exact scanning (FES) by \cite{ma2016fisher} constructed multi-scale scan statistics that are particularly effective at detecting local dependency through Fisher's exact tests over rectangle scanning windows.

Most of the above nonparametric tests enjoy the property of universal consistency against any particular form of dependence. Formally, this universality means that for any specific copula distribution $\Pb_1 \neq \Pb_0$, as $n \rightarrow \infty$, the test for the problem $ H_0: \Pb_{(U,V)} = \Pb_0 ~~v.s.~~ H_1: \Pb_{(U,V)}=\Pb_1$ has an asymptotic power of $1$. However, one important problem in many distribution-free tests is the lack of uniformity. To see this, we consider the total variation (TV) distance $TV(\cdot,\cdot),$ which is defined by $TV(\Pb,\Qb)=\sup_{S \in \cF}|\Pb(S)-\Qb(S)|,$ where $\cF$ is a $\sigma$-algebra of the sample space. The uniform consistency of nonparametric dependence detection w.r.t. the TV distance is to be consistent for any alternative which is certain distant from independence, i.e.,
{\small
\begin{equation}\label{eq:test_ind_unif}
H_0: \Pb_{(U,V)} = \Pb_0 ~~v.s.~~ H_1: TV(\Pb_{(U,V)},\Pb_0) \ge \delta
\end{equation}
}for some $0<\delta\le 1.$ For the testing problem in \eqref{eq:test_ind_unif}, although many tests are universally consistent, we show in Section~\ref{sec: sudden} and Theorem~\ref{thm: sudden} the non-existence of a test that is uniformly consistent w.r.t. the TV distance. The uniformity issue is due to the fact that the space of $H_1$ is large. Said another way: When two variables are not independent, there are so many ways they can be dependent. In practice, having this non-uniform consistency problem means having ``blind spots'' in dependence detection for a given sample size, i.e., having very low power for many forms of dependency, especially nonlinear ones. Note that nonlinear forms of dependence are ubiquitous in sciences, for example laws in physics defined by differential equations. Therefore, avoiding the power loss due to the non-uniform consistency problem in nonparametric dependence detection means having robust power against a large class of alternatives and improving the ability of discovering novel relationships in many areas of science.

Because of the impossibility of testing \eqref{eq:test_ind_unif} with uniform consistency w.r.t. the TV distance (Theorem~\ref{thm: sudden}), to avoid such power loss, we propose to test approximate independence through a \textit{filtration approach}. Such a filtration is constructed by the $\sigma$-fields generated by binary variables from marginal binary expansions which jointly approximate the copula distribution. Similar filtration ideas are nicely described in \cite{liu2014simpson,liu2016individual} in studying the Simpson's Paradox. The approximation idea is also related to the ``probably approximately correct'' (PAC) approach in machine learning \citep{valiant1984theory}. We explain the details in Section~\ref{subsec: be}.

We note here that although many other ways of filtration approximations are available, there are a few important advantages of the proposed binary expansion filtration that facilitate studies of dependence.

(a) The $\sigma$-field generated by binary variables is \emph{finite}.

(b) Two binary variables are independent if and only if they are \emph{uncorrelated}.

We call the statistics that are functions of the Bernoulli variables from the above filtration approximation binary expansion statistics (BEStat), and we call the testing framework on the corresponding approximate independence the binary expansion testing (BET) framework. This approach leads to studies of contingency tables from discretizations. Although classical tests such as the $\chi^2$ tests \citep{lehmann2006testing} are readily available, they have some drawbacks: (a) the exponentially growing degrees of freedom that would affect the power, and (b) the unclear interpretability of dependence when the independence hypothesis is rejected. To improve on these two issues, we consider reparametrization of the likelihood of the contingency tables through a novel binary interaction design (BID) equation (Theorem~\ref{thm: bid}), which connects the study of dependence to the Hadamard transform in signal processing. Through this connection, the interactions of binary variables in the filtration are shown to be complete sufficient statistics for dependence. By utilizing these interactions, we convert the dependence detection problem to a multiple testing problem. Statistically speaking, the benefits of the above approach are summarized below:
\vspace*{-0.2\baselineskip}
\begin{enumerate}[(a)]
\item The Hadamard transform provides new insights for the analysis of any contingency table whose size is a power of 2. Compared to the conventional parametrization, the novel parameters marginal interaction odds ratios (MIOR) and cross interaction odds ratios (CIOR) separate the marginal and joint information, and CIORs being 1 is equivalent to independence. As an analogy, the CIORs are to contingency tables as the correlations are to multivariate normal distributions. See Theorem~\ref{thm: be_log_ior} and Theorem~\ref{thm: be_ior}.
\item The symmetry statistics from the reparametrization are shown to be complete sufficient statistics for dependence. They are identically distributed and are uncorrelated under the null. See Theorem~\ref{thm: symmetry}, Theorem~\ref{thm: symmetry_sample} and Theorem~\ref{thm: ortho}.
\item As a consequence of the above properties, the multiple testing procedure is shown to be minimax in the sample size requirement for reliable power. See Theorem~\ref{thm: optimality}.
\item Upon rejection of independence, the largest absolute symmetry statistic and the corresponding cross interaction provide clear interpretation of the dependency.
\end{enumerate}

\noindent Although theories for copula and contingency tables are well-developed, we are not aware of similar approach or results in statistical literature.

We also note that the BEStat approach is closely related to computing. In current computing systems, each decimal number is coded as a sequence of binary bits, which is exactly the binary expansion of that number. This connection means that one can carry out the BEStat procedures by operating directly over bits. Since bitwise operations are one of the most efficient operations in current computing systems, we are able to develop computationally efficient implementations of the proposed method. The detailed algorithm is described in a separate paper \citep{zhao2019fast}, and it improves the speed of existing methods by orders of magnitude.

This paper is organized in as follows: Section~\ref{sec: sudden} explains the problem of non-uniform consistency. Section~\ref{sec: bestat} introduces the concept and basic theory in the framework of BEStat and BET. Section~\ref{sec: maxbet} studies the Max BET procedure and its properties. Section~\ref{sec: computing} connects the BEStat framework to current computing system. Section~\ref{sec: numerical}, Section~\ref{sec: stars} and Section~\ref{sec: tcga} illustrate the procedure with simulated and real data studies. Section~\ref{sec: summary} concludes the paper with discussions of future work. The proofs can be found in the supplementary file.

\section{Motivation: Non-Uniform Consistency}\label{sec: sudden}
To explain the problem of non-uniform consistency, we develop the following example of the bisection expanding cross (BEX). Many existing methods suffer substantial power loss under this example due to this problem, which can be avoided through the binary expansion statistics proposed in Section~\ref{sec: bestat} and Section~\ref{sec: maxbet}.

\begin{figure}[hhhh]
\begin{center}
\includegraphics[width=\textwidth]{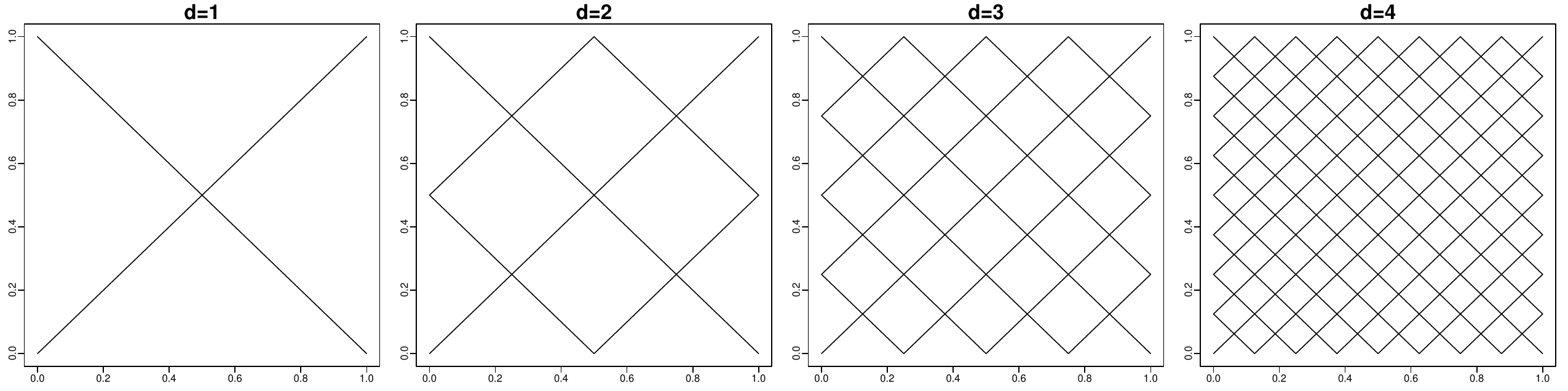}
\end{center}
\caption{
The bisection expanding cross (BEX) at level $d=1,\ldots,4$.}\label{fig: bex}
\end{figure}

We call the following sequence of one-dimensional manifolds in $[0,1]^2$ the bisection expanding cross (BEX). These manifolds can be defined through the implicit function $\gamma_d(x,y)=0$ for every integer $d>0$: $BEX_d = \{(x,y) \in [0,1]^2: \gamma_d(x,y)=0\}$, where
{\footnotesize
\begin{equation*}
\gamma_d(x,y) = \sum_{i=1}^{2^{d-1}}\sum_{j=1}^{2^{d-1}} \bigg(\bigg|x-{i \over 2^{d-1}}+{1 \over 2^d}\bigg|-\bigg|y-{j \over 2^{d-1}}+{1 \over 2^d}\bigg|\bigg){\rm I}\bigg(\bigg|x-{i \over 2^{d-1}}+{1 \over 2^d}\bigg| \le {1 \over 2^d}\bigg){\rm I}\bigg(\bigg|y-{j \over 2^{d-1}}+{1 \over 2^d}\bigg| \le {1 \over 2^d}\bigg).
\end{equation*}
}The BEX structure is illustrated in Figure~\ref{fig: bex}, where the first four levels are plotted. Graphically, this grid can be regarded as a space-filling fractal by recursively expanding the bisector of the four ``arms'' of $BEX_1$ until intersection.

Now we consider the random variables $(X_d, Y_d)$ that are uniformly distributed over $BEX_d$ whose joint distribution is denoted by $\Pb_d$. The properties of these distributions are summarized in the following proposition.
\begin{proposition}\label{prop: unif_bex}
\ \begin{enumerate}[(a)]
    \item $X_d$ and $Y_d$ are marginally $Uniform[0,1]$ for any $d$.
    \item $\gamma_d(X_d,Y_d)=0$ for any $d$, i.e., the joint distribution of $(X,Y)$ is degenerate. In particular, $TV(\Pb_d,\Pb_0)=1$ for any $d$.
    \item $\forall (x,y) \in [0,1]^2$, as $d \rightarrow \infty$, $|\Pb_d(X_d \le x,Y_d\le y) - \Pb_d(X_d \le x)\Pb_d(Y_d \le y)|\rightarrow 0$.
  \end{enumerate}
\end{proposition}
Part (b) and part (c) of Proposition~\ref{prop: unif_bex} seem to contradict each other: Part (b) says that the joint distribution of $X_d$ and $Y_d$ is far away from independence in the TV distance, thus they are strongly non-independent. Yet, part (c) claims that when $d$ is large, $X_d$ and $Y_d$ are nearly independent. Indeed, the BEX shows that despite a TV distance of 1, degenerate distributions can be arbitrarily close to independence. We shall explain this paradox in Section~\ref{subsec: interpretation}. This paradox also lead to a challenge to testing methods: Given a finite sample, can we effectively distinguish any form of dependency from independence?

Unfortunately, for any testing method, the answer is negative. Intuitively speaking, this is because for any given test with a given samples size $n$, one can keep expanding the BEX until it is so close to independence that this test becomes powerless. This example thus illustrates the problem of non-uniform consistency of the test in \eqref{eq:test_ind_unif}: No test can be uniformly consistent against all forms of dependence, not even all levels of the BEX, for which $\delta=1$ in \eqref{eq:test_ind_unif}. See Theorem~\ref{thm: sudden} below.

The power loss due to non-uniform consistency can be severe. For example, simulations (see Section 1.1 in the supplementary file) show that many CDF based and kernel based tests are powerless in detecting BEX at level $4$ even when the sample size is as high as $20000$. Note that with such a large sample, the BEX structure and the dependency can be clearly observed in the scatterplot by naked eyes. However, many existing tests cannot distinguish it from independence.

We make a few remarks about the BEX example before proceeding.

(a) The BEX is closely related to many research problems such as the chessboard detection in computer vision \citep{forsyth2002computer}.

(b) The BEX is not the first example that a sequence of degenerate distributions converges to independence. The earliest example we could find is in \cite{kimeldorf1978}. There are also other interesting and useful fractal applications in statistics such as \cite{craiu2005multiprocess,craiu2006meeting}. The basis of the BEX example is a classical result in \cite{Vitale1990}. We construct the BEX paradox due to its fractal structure which explains the problem of non-uniform consistency.

(c) The non-uniform consistency shown with the BEX is specifically for our choice of the TV distance between distributions. There are many other distances \citep{tsybakov2008introduction}, and a different choice of distance could lead to a different test statistic and different results on uniform consistency. We choose the TV distance because (1) it is a widely used distance in literature, (2) it is equivalent to many other distances, and (3) it is convenient for the analysis in our binary expansion approach. Therefore, throughout this paper, we focus on the TV distance, and all results about uniform consistency are w.r.t. the TV distance. In particular, we provide a formal statement of the problem of non-uniform consistency w.r.t. the TV distance below:

\begin{theorem}\label{thm: sudden}
Consider the testing problem in \eqref{eq:test_ind_unif}. For any finite number of i.i.d. observations $n$, for any test that has a Lebesgue measurable critical region $C_n \subset \RR^{2n}$ with $\Pb_{H_0}(\partial C_n)=0$ and $\Pb_{H_0}(C_n) \le \alpha$, $\forall \epsilon >0$, there exists a bivariate distribution $F_n \in H_1$ and $\Pb_{F_n}(C_n) \le \alpha +\epsilon$.
\end{theorem}

The message of Theorem~\ref{thm: sudden} is that in a distribution-free setting without any assumption on the joint distribution, dependence is not a tractable target. The intractability comes from the fact that without a model of the joint distribution, there is no parameter to characterize and identify the underlying form of dependency. Therefore, there is no target for inference about dependence from a test or any other statistical method. Although one can develop good measures of dependence such as distance correlation, GMC, HSIC and MIC, etc., such measures cannot make the joint distribution identifiable. Therefore, they can never replace the role of parameters in statistical inference about dependence. This fact motivates the following three key elements in the BEStat approach and the BET framework:
\begin{enumerate}[(a)]
\item Rather than one test of independence, we will study dependence through a carefully designed sequence of tests based on a filtration to achieve \textit{universality}.
\item For every test statistic in the sequence, there is an explicit well-defined set of parameters as the target for inference to achieve \textit{identifiability}.
\item At every step in the sequence, the test is consistent against all alternatives which are $\delta$-away from independence in the TV distance to achieve \textit{uniformity}.
\end{enumerate}
The above BET framework can help explain the seeming paradox in the BEX example, and the proposed test can have high power against this dependency. See Section~\ref{subsec: interpretation}.

\section{The Basic Theory of Binary Expansion Statistics}\label{sec: bestat}
\subsection{Binary Expansion Filtration}\label{subsec: be}
The considerations in Section~\ref{sec: sudden} necessitate a multi-scale binning approach to study dependence. For the dependence detection problem, this multi-scale approach means to test some approximate independence rather than the exact hypothesis in \eqref{eq:test_ind_unif}. We study the known marginal CDF case first, for which we develop such a multi-scale framework through the following classical result on the binary expansion of a uniform random variable \citep{kac1959statistical}:
\begin{theorem}\label{thm: be}
If $U \sim Uniform[0,1],$ then $U=\sum_{k=1}^\infty { A_k \over 2^k}$ where $A_k \stackrel{i.i.d.}{\sim} Bernoulli(1/2).$
\end{theorem}

The binary expansion in Theorem~\ref{thm: be} decomposes the information about $U$ into information from independent Bernoulli $A_k$'s. $A_k$'s can be also regarded as indicator functions of $U$. For example, $A_1={\rm I}(U\in (1/2,1])$, $A_2={\rm I}(U \in (1/4,1/2]\cup (3/4,1])$, see \cite{kac1959statistical}. To study the dependence between $U$ and $V$, we consider the binary expansion of both $U$ and $V$: $U=\sum_{k=1}^\infty { A_k \over 2^k}$ and $V=\sum_{k=1}^\infty { B_k \over 2^k}$ where $A_k\stackrel{i.i.d.}{\sim} Bernoulli(1/2)$ and $B_k\stackrel{i.i.d.}{\sim} Bernoulli(1/2)$.

Note that if we truncate the binary expansions of $U$ and $V$ at some finite depths $d_1$ and $d_2$ respectively, $U_{d_1}=\sum_{k=1}^{d_1} { A_k \over 2^k}$ and $V_{d_2}=\sum_{k=1}^{d_2} { B_k \over 2^k}$, then $U_{d_1}$ and $V_{d_2}$ are two discrete variables that can take $2^{d_1}$ and $2^{d_2}$ possible values respectively. Moreover, as $d_1, d_2\rightarrow \infty$, $|U_{d_1}-U|=O_p(2^{-d_1})$ and $|V_{d_2}-V|=O_p(2^{-d_2})$. In particular,
{\small
\begin{equation}\label{eq: beuvd_joint}
\|(U_{d_1},V_{d_2})- (U,V)\|_2=O_p(2^{-\min\{d_1,d_2\}}).
\end{equation}
}
\vspace*{-\baselineskip}

The above considerations are apparent if one regards the truncations as a filtration generated by $\{A_k\}_{k=1}^{d_1}$ and $\{B_k\}_{k=1}^{d_2}$ for each $d_1, d_2\ge 1$. Indeed, the filtration idea is a consequence of George Box's aphorism ``All models are wrong, but some are useful.'' At every $d_1$ and $d_2$, the probability model of $(U_{d_1},V_{d_2})$ is a ``wrong'' model for the joint distribution $(U,V)$. However, the ``wrong'' model of $(U_{d_1},V_{d_2})$ can be very useful in many ways. In particular, we show below how the three key elements described at the end of Section~\ref{sec: sudden} are achieved from this approach:

(a) \textit{Universality}: The important message from \eqref{eq: beuvd_joint} is that one can approximate the joint distribution of and hence the dependence in $(U,V)$ through that in $(U_{d_1},V_{d_2}).$ Although the dependence in the joint distribution of $(U,V)$ can be arbitrarily complicated, when $d_1$ and $d_2$ are large, we expect a good approximation from discrete variables $(U_{d_1}, V_{d_2})$ where the approximation error is exponentially small. In terms of testing independence, this means although the joint distribution of $(U,V)$ can be arbitrarily close to independence, due to the filtration feature of the sequence, one can always detect the dependence when $d_1$ and $d_2$ are large to achieve universality.

(b) \textit{Identifiability}: As we explained in Section~\ref{sec: sudden}, one crucial challenge in distribution-free dependence detection is identifiability. Without models and parameters, dependence is not a tractable target. On the other hand, $(U_{d_1},V_{d_2})$ can only take a finite $2^{d_1+d_2}$ possible values, which leads to a partition of the scatterplot of data into a $2^{d_1} \times 2^{d_2}$ contingency table. With this consideration, the truncation of the binary expansions turns the problem on dependence, which is unidentifiable under the distribution-free setting, into a problem over a contingency table, which is fully identifiable. In terms of testing, when we begin without any assumptions about the joint distribution, there is no explicit way to write out the alternative likelihood under dependence. However, at each depths $d_1$ and $d_2$, due to the discreteness, the class of alternative distributions is restricted to those over the contingency table, which has an explicit distribution and has cell probabilities as identifiable parameters for inference \citep{agresti2011categorical,fienberg2007analysis}.

(c) \textit{Uniformity}: As a consequence of identifiability, we can avoid the problem of non-uniform consistency described in Section~\ref{sec: sudden}. At any depths $d_1$ and $d_2$, one can write out the TV distance between an alternative distribution and the null distribution in terms of the cell probabilities in the contingency table model. We are thus able to show the consistency and optimality of the proposed Max BET procedure in Section~\ref{subsec: optimality} for alternative distributions whose TV distances from the independence null is at least $\delta$, for any $\delta>0$.

The above considerations motivate us to propose the binary expansion statistics in studying the dependence between $U$ and $V$ in a distribution-free setting. Formally, we define binary expansion statistics as follows:
\begin{definition}\label{def: bestat}
  We call statistics as functions of finitely many Bernoulli variables from marginal binary expansions the binary expansion statistics (BEStat).
\end{definition}
Similarly, for the problem of detecting dependence from independence in a distribution-free setting, we define the binary expansion testing framework as follows.
\begin{definition}\label{def: bet}
  We call the testing framework based on the binary expansion filtration approximation up to certain depth the binary expansion testing (BET).
\end{definition}
In the context of testing independence in bivariate distributions, the BET at depths $d_1$ and $d_2$ is to test the independence of $U_{d_1}$ and $V_{d_2}$, which we refer to as $(d_1,d_2)$-independence and which is equivalently defined in \cite{ma2016fisher} for scanning statistics. Formally, denote the bivariate uniform distribution over $\{{0 \over 2^{d_1}},\ldots,{2^{d_1}-1 \over 2^{d_1}}\} \times \{{0 \over 2^{d_2}},\ldots,{2^{d_2}-1 \over 2^{d_2}}\}$  by $\Pb_{0,d_1,d_2}$. For some $0<\delta\le 1, $ we consider
{\small
\begin{equation}\label{eq:test_be_ind}
H_{0,d_1,d_2}: \Pb_{(U_{d_1},V_{d_2})} = \Pb_{0,d_1,d_2} ~~v.s.~~ H_{1,d_1,d_2}: TV(\Pb_{(U_{d_1},V_{d_2})},\Pb_{0,d_1,d_2}) \ge \delta.
\end{equation}
}
\vspace*{-\baselineskip}

Not rejecting the null hypothesis in the BET at depths $(d_1,d_2)$ thus indicates that there is no strong evidence against the null hypothesis of independence between $U$ and $V$ up to depths $d_1$ and $d_2$ in the binary expansions. Note that this interpretation is weaker than claiming independence between $U$ and $V$: The dependence can occur at some larger $(d_1,d_2)$ in the $O_p(2^{-\min{\{d_1,d_2\}}})$ remainder term in \eqref{eq: beuvd_joint}. However, as described in Section~\ref{sec: sudden}, claiming exact independence with finite samples and without any restriction on the alternative is impossible. On the other hand, this weaker hypothesis of approximate independence helps us to avoid the uniform consistency problem in the dependence detection under the distribution-free setting and provides reliable power for a large class of alternatives. To see the gains from this trade-off, one can compare our results in Section~\ref{subsec: optimality} with those in Section~\ref{sec: sudden}.

We remark here that the filtration in approximating dependence is not unique. For example, one can consider the filtration corresponding to orthogonal polynomials rather than the binary expansion. However, the $\sigma$-field in the binary expansion filtration has a few important advantages to facilitate studies of dependence.

(a) Finiteness of $\sigma$-fields: For the $\sigma$-field at each depths $d_1$ and $d_2$, the number of events is $2^{d_1+d_2}-1,$ which is finite. This is because interactions of binary variables are at most binary. If we consider some other filtration (for example orthogonal polynomials) for the approximation of dependence, then the $\sigma$-field might not be of finitely many events and can be much more complicated.

(b) Uncorrelatedness implying independence: Although uncorrelatedness usually does not imply independence, it is well known that it does for two binary variables. This property can greatly simplify studies of dependence in filtration. Again, if we consider some other filtration (for example orthogonal polynomials) for the approximation of dependence, then quantifying the dependence between variables in the $\sigma$-field can be much more complicated.

The above considerations also work similarly for the case when the marginal distributions are unknown. To study the binary expansion in this case, suppose the sample size is $n=2^K$ for some $K>0$ for easy explanation. With the marginal empirical CDF transformations, the $i$-th observation in the empirical copula are $\hat{U}_i$ and $\hat{V}_i$ whose marginal distribution is $Uniform\{{1 \over 2^K},\ldots,{2^K \over 2^K}\}.$ Now let $\hat{A}_{1,i}= {\rm I}(\hat{U}_i\in (1/2,1]), \ldots, \hat{A}_{K,i}={\rm I}(\hat{U}_i \in \cup_{k'=1}^{2^{K-1}}({2k'-1\over 2^K}, {2k' \over 2^K}] ).$ It is easy to see that for each fixed $i$, $\hat{A}_{k,i}$'s are independent, and $\hat{U}_i= {1 \over 2^K}+\sum_{k=1}^K { \hat{A}_{k,i} \over 2^k}.$ Therefore, the binary expansion filtration can be similarly defined, and the BET at depths $d_1$ and $d_2$ is to test the independence of $\hat{U}_{d_1,i}= \sum_{k=1}^{d_1} { \hat{A}_{k,i} \over 2^k}$ and $\hat{V}_{d_2,i}=\sum_{k=1}^{d_2} { \hat{B}_{k,i} \over 2^k}$:
{\small
\begin{equation}\label{eq:test_be_ind_sample}
H_{0,d_1,d_2}: \text{For each }i, \hat{U}_{d_1,i}\text{ and }\hat{V}_{d_2,i}\text{ are independent}.
\end{equation}
}The interpretation of this null hypothesis is that for each observation, the row assignment and column assignment to the contingency table are independent, as in classical categorical data analysis \citep{agresti2011categorical,fienberg2007analysis}. When $\hat{U}_{K,i}$ and $\hat{V}_{K,i}$ are independent for each $i$, the observed ranks are independent.

We explain the details of these tests in Section~\ref{subsec: likelihood} and Section~\ref{sec: maxbet}. We remark here that although copula theory is well developed \citep{nelsen2007introduction}, we are not aware of any filtration approach in the literature. We also remark here that tests of approximate independence are also considered in a very recent paper \citep{ma2016fisher} for scanning purposes, in which a filtration idea is implicitly described. In this paper, our goal is to formally develop the framework of binary expansion statistics. We shall compare the theory and methods in both papers in Section~\ref{subsec: other}.

\subsection{Revisiting the Classical Theory for Contingency Tables}\label{subsec: likelihood}
We start our analysis by first revisiting the model and theory of a general contingency table with $r$ rows and $c$ columns of $n$ i.i.d. samples. The parameters of interest are $\bp=\{p_{ij}, i=1,\ldots,r, j=1,\ldots,c\}$, and the cell counts are $\bn=\{n_{ij}\}$. The only constraint is on the totals $\sum_{i,j} p_{ij}=1$ and $\sum_{i,j}n_{ij} = n.$ Two most important models for the likelihood are as follows \citep{agresti2011categorical,fienberg2007analysis}:

(a) When there is no restriction on marginal totals, the joint distribution of the cell count vector $\bN$ is multinomial (with the convention $0^0=1$): With $C_1(\bn)= {n! \over \prod_{i,j} n_{ij}!},$
{\small
\begin{equation}\label{eq: likelihood_multinomial}
p(\bN=\bn|\bp)= C_1(\bn)\prod_{i,j} p_{ij}^{n_{ij}}.
\end{equation}
}
\vspace*{-\baselineskip}

(b) Condition on positive row and column totals $\bn_r=\{n_{i\cdot}=\sum_{j} n_{ij},i=1,\ldots,r \}$ and $\bn_c=\{n_{\cdot j}=\sum_{i} n_{ij},j=1,\ldots,c \}$, for $i<r$ and $j<c$, with the reparametrization $\theta_{ij}={p_{ij}p_{rc} \over p_{ic}p_{rj}}$ and normalizing constant $h_1(\bn_r,\bn_c,\btheta)$, we have $p(\bN=\bn|\btheta,\bn_r,\bn_c)= C_1(\bn)h_1(\bn_r,\bn_c,\btheta) \prod_{i,j}\theta_{ij}^{n_{ij}}$ \citep{cornfield1956statistical}. Note that under independence $\theta_{ij}=1$, and the distribution is (central) multivariate hypergeometric
{\small
\begin{equation}\label{eq: cmvhypergeo}
p(\bN=\bn|\bn_r,\bn_c)=  C_1(\bn)h_1(\bn_r,\bn_c)={\prod_i n_{i \cdot}! \prod_j n_{\cdot j}! \over n! \prod_{i,j} n_{ij}!},
\end{equation}
}With the above distributions, tests of independence for a contingency table can be done through classical methods such as $\chi^2$ tests, Fisher's exact tests, and likelihood ratio tests (LRT). For the nonparametric dependence detection problem, the BET with these tests are uniformly consistent for any depths $d_1$ and $d_2$. However, these classical methods have two important limitations on power and interpretability:

(a) The minimal sample size for classical tests to have reliable power is known \citep{agresti2011categorical,fienberg2007analysis} to be about the size of the contingency table $O(2^{d_1+d_2}).$ However, recent developments \citep{acharya2015optimal} show that the optimal lower bound of this sample size requirement is $O(2^{d_1+d_2 \over 2})$. This result indicates that classical tests may suffer substantial power loss in dependence detection, especially when $d_1$ and $d_2$ are large. For a well-known example, when the contingency table contain many empty cells, LRT and $\chi^2$ tests will fail to work.

(b) The rejections from classical tests are not very interpretable. Even if we can claim significant dependence with a classical test, the test does not provide information about how the variables are dependent.

One intuition of the above limitations in classical tests is that each cell in a contingency table is considered in an isolated manner, thus the information between cells is somehow lost. To improve classical tests, we consider grouping the cells together to improve the power and interpretability. Such grouping process is effectively achieved through the binary interaction design described in Section~\ref{subsec: interaction}.

\subsection{Binary Interaction Design: Reparametrization of the $2^{d_1} \times 2^{d_2}$ Contingency Table Likelihood}\label{subsec: interaction}
We now turn to the case when the contingency table is generated by the binary expansion up to depths $d_1$ and $d_2$ as described in Section~\ref{subsec: be}, so that the table has $2^{d_2}$ rows and $2^{d_1}$ columns (assuming $U$ on the horizontal axis and $V$ on the vertical axis). To provide a general theory for contingency tables, in this subsection we \emph{do not} restrict the total probability of each row and column being the same (which happens when $A_i$'s and $B_j$'s are both i.i.d. $Bernoulli(1/2)$). However, in this subsection, we shall assume that all cell probabilities are positive.

To combine the cell information, we consider the $\sigma$-field generated from the binary expansion filtration. We explain in the known marginal distribution case first since it is similar for the unknown marginal distribution case. With $d_1$ Bernoulli variables $A_k,k=1,\ldots,d_1$ and another $d_2$ Bernoulli variables $B_k,k=1,\ldots,d_2$ (again in this subsection we \emph{do not} assume them to be independent and symmetric), consider two general discrete variables defined by $U_{d_1}=\sum_{k=1}^{d_1} { A_k \over 2^k}$ and $V_{d_2}=\sum_{k=1}^{d_2} { B_k \over 2^k}.$ The $\sigma$-field here is $\sigma(U_{d_1},V_{d_2})=\sigma(A_1,\ldots,A_{d_1},B_1,\ldots,B_{d_2})$ and is generated by $2^{d_1+d_2}-1$ Bernoulli variables resulting from interactions between $A_i$'s and $B_j$'s. We shall use the equivalent binary variables $\dot{A}_i=2A_i-1$ and $\dot{B}_j=2B_j-1$ since the interaction between them can be conveniently written as products. For example, the event $\{A_1=1,B_1=1\}\cup \{A_1=0,B_1=0\}$ is equivalent to the event $\{\dot{A}_1\dot{B}_1=1\}$.

\begin{figure}[hhhh]
\begin{center}
\includegraphics[width=\textwidth]{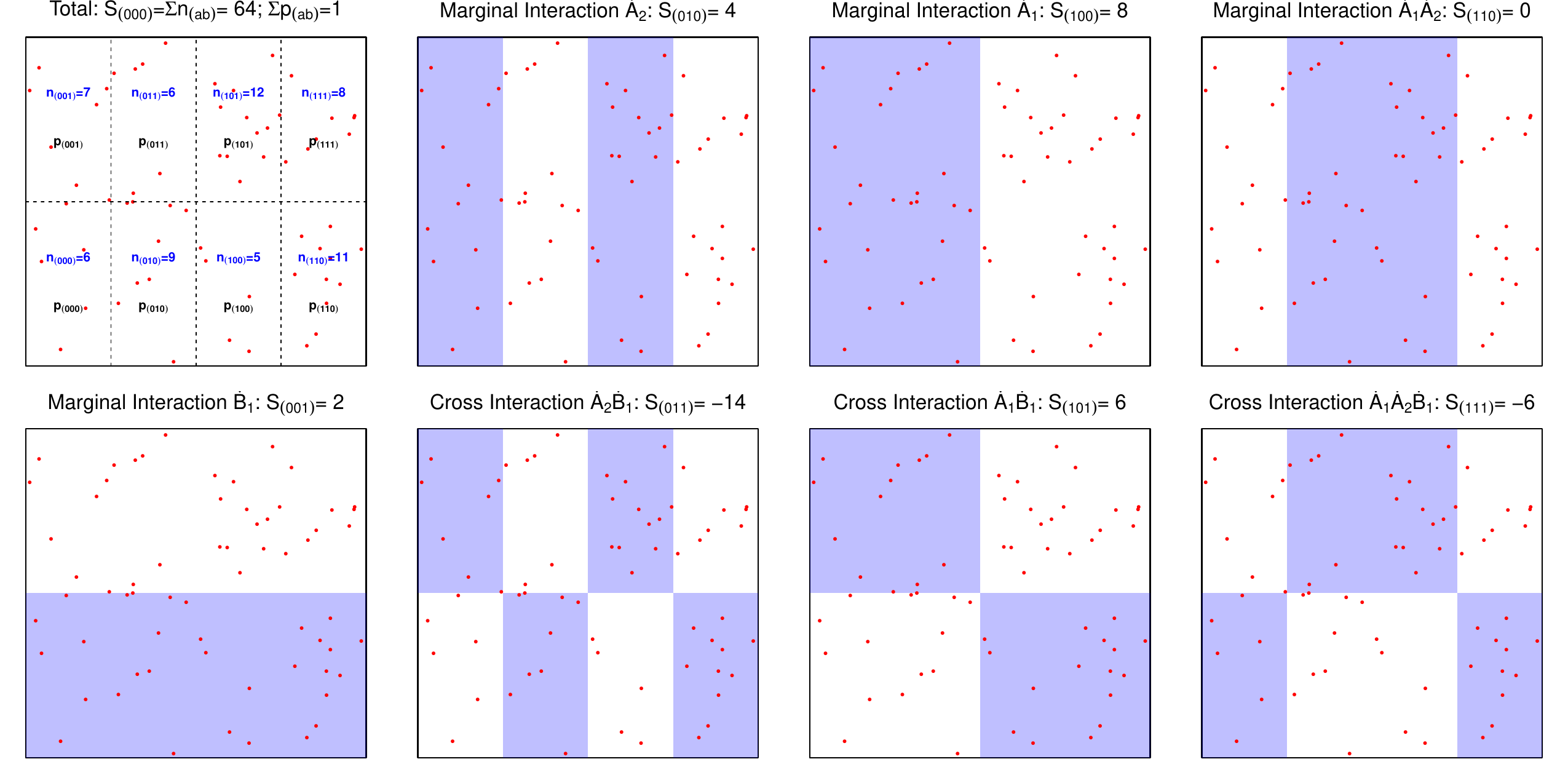}
\end{center}
\caption{
The binary interaction design (BID) at depths $d_1=2$ and $d_2=1$ with $n=64$ observations. The number of observations in each cell is presented in the top left plot. There are 7 non-trivial binary variables in the $\sigma$-field, whose positive regions are in white and whose negative regions are in blue. Symmetry statistics $S_{(\ba\bb)}$ are calculated for these 4 marginal interactions and 3 cross interactions. For example, $S_{(011)}=n_{(111)}-n_{(110)}-n_{(101)}+n_{(100)}+n_{(011)}-n_{(010)}-n_{(001)}+n_{(000)}=-14$.}\label{fig: bid}
\end{figure}

Note that each of these binary interaction variables leads to a partition of the unit square $[0,1]^2$ and two groups of cells according to whether the interaction is positive. Moreover, for each interaction in the $\sigma$-field, the number of cells in the regions where it takes value $1$ (and $-1$) is exactly $2^{d_1+d_2-1}.$ This fact can be explained by the BID equation (Theorem~\ref{thm: bid}) below, and it facilitates the definition of interaction odds ratio (IOR) as in Definition~\ref{def: ior} as well as the reparametrization with IOR. The IORs group the cell information together and separate the marginal and joint information in the multinomial likelihood. See Figure~\ref{fig: bid}.

Note also that the $2^{d_1+d_2}-1$ binary variables in the $\sigma$-field can be categorized into two classes: The variables of the form $\dot{A}_{k_1}\ldots \dot{A}_{k_r}$ or $\dot{B}_{k'_1}\ldots \dot{B}_{k'_t}$ will be referred to as \emph{marginal interactions} since they only involve the marginal distributions. On the other hand, the variables of the form $\dot{A}_{k_1}\ldots \dot{A}_{k_r}\dot{B}_{k'_1}\ldots \dot{B}_{k'_t}$ with $r,t>0$ will be referred to as \emph{cross interactions} since they contain information of both $U_{d_1}$ and $V_{d_2}$.

In explanation of the theory, we use the following binary integer indexing for related quantities: Denote the Bernoulli random vectors in the binary expansion by $\bA=(A_{1},\ldots,A_{d_1})$ and $\bB=(B_{1},\ldots,B_{d_2})$, and denote vectors of length $d_1$ and $d_2$ with entries $0$'s and $1$'s by $\ba$ and $\bb.$ The probability of each of the $2^{d_1+d_2}$ cells can then be written as $p_{(\ba\bb)} = \Pb(\bA=\ba, \bB=\bb)$ with $(\ba\bb)$ being the concatenation of $\ba$ and $\bb.$ Now let the integer $c$ determined by $c=\sum_{i=1}^{d_1} a_{i} 2^{d_1+d_2-i}+\sum_{j=1}^{d_2} b_j 2^{d_2-j}.$ Let $\bp$ be the $2^{d_1+d_2}$-dimensional vector of probabilities whose $(2^{d_1+d_2}-c)$-th entry is $p_{(\ba\bb)}.$

For the binary variables in $\sigma(\dot{A}_1,\ldots,\dot{A}_{d_1},\dot{B}_1,\ldots,\dot{B}_{d_2}),$ we also denote their expected values with binary integer index as follows. For $\Eb[\dot{A}_{k_1}\ldots \dot{A}_{k_r}\dot{B}_{k'_1}\ldots \dot{B}_{k'_t}], r=1,\ldots,d_1, t=1,\ldots,d_2$, we denote it by $E_{(\ba\bb)}$ where $\ba$ is a $d_1$-dimensional binary vector with 1's at $k_1,\ldots,k_r$ and are 0's otherwise, and $\bb$ is a $d_2$-dimensional binary vector with 1's at $k'_1,\ldots,k'_t$ and are 0's otherwise. Note here that $E_{(\zero \zero)}=\Eb[1]=1.$ We also write the interaction as a product of binary variables $\dot{A}_{k_1}\ldots \dot{A}_{k_r}\dot{B}_{k'_1}\ldots \dot{B}_{k'_t}$ as $\dot{A}_{\ba}\dot{B}_{\bb}.$ With $c$ defined in the previous paragraph, let $\bE$ be the $2^{d_1+d_2}$-dimensional vector of expected values whose $(c+1)$-th entry is $E_{(\ba\bb)}.$

The above notation also applies to observed quantities: With the total $n$ observations, the cell counts are denoted by $n_{(\ba\bb)}.$ The collection of all $n_{(\ba\bb)}$'s is denoted by $\bN$ and is indexed as in $\bp$. We also denote the sum of observed binary interaction variables by $S_{(\ba\bb)}=\sum_{i=1}^n \dot{A}_{\ba,i} \dot{B}_{\bb,i}$ with $S_{(\zero \zero)}=n.$ The collection of all $S_{(\ba\bb)}$'s is denoted by $\bS$ and is indexed as in $\bE$. We shall refer $S_{(\ba\bb)}$ as \emph{the symmetry statistic} for $\dot{A}_{\ba}\dot{B}_{\bb}$ as they can be regarded as the differences between the numbers of points in positive and negative regions. Thus, $S_{(\ba\bb)}$ is a statistic about symmetry. See Figure~\ref{fig: bid}.

With the above notation, we establish the equation connecting the contingency table distribution and the interactions of binary variables in the $\sigma$-field. The equation is established through $\Hb=\Hb_{2^{d_1+d_2}}$ being the Sylvester's construction of Hadamard matrix \citep{sylvester1867lx}. We shall refer this equation as the binary interaction design (BID) equation (name coined in \cite{zhao2019fast}).
\begin{theorem}\label{thm: bid}
\ \begin{enumerate}[(a)]
\item Population version of the BID equation: $\bE=\Hb \bp.$
\item Sample version of the BID equation: $\bS=\Hb \bN.$
\end{enumerate}
\end{theorem}
The Hadamard matrix $\Hb$ is referred to as Walsh matrix in literature of signal processing, where a linear transformation with $\Hb$ as in Theorem~\ref{thm: bid} is referred to as the Hadamard transform \citep{lynn1973introduction,golubov2012walsh,harmuth2013transmission}. The earliest referral to the Hadamard matrix we found in statistical literature is \cite{Pearl1971}. The Hadamard matrix is also closely related to the orthogonal full factorial design \citep{box2005statistics,cox2000theory}. In the context of dependence detection, this transform maps the cell domain (in $\bp$ or $\bN$) to the interaction domain (in $\bE$ or $\bS$). Thus, the information in individual cells can be grouped together to provide information about global dependency. Although theory and methods for contingency tables are well-developed, we are not aware of similar approach in related literature.

To see the importance of the BID equation and the symmetry statistic $S_{(\ba\bb)}$, we introduce some more notation here. We label the first to $2^{d_1+d_2}$-th row (and column) of $\Hb$ with binary integer indices from $(\zero_{d_1+d_2})$ to $(\one_{d_1+d_2})$. Denote $\overline{(\ba\bb)}=(\one\one)-(\ba\bb)$ to be the \emph{binary conjugate}, or logical negation of $(\ba\bb)$, i.e., $\overline{(010)}=(101).$ With the above notation, we summarize some useful properties of the Hadamard matrix $\Hb_{2^{d_1+d_2}}$ in the following proposition \citep{golubov2012walsh}.
\begin{proposition}\label{prop: hadamard}
  \ \begin{enumerate}[(a)]
  \item $\Hb_{2^{d_1+d_2}}$ is symmetric. The entry in $\Hb_{2^{d_1+d_2}}$ at the $(\ba'\bb')$-th row and $(\ba\bb)$-th column is $(-1)^{(\ba'\bb')^T(\ba\bb)}.$
  \item $\Hb_{2^{d_1+d_2}}$ has orthogonal columns: $\Hb_{2^{d_1+d_2}}^{-1} = {1 \over 2^{d_1+d_2}}\Hb_{2^{d_1+d_2}}.$
  \item Hadamard matrices can be defined recursively: $\Hb_{2^{d_1+d_2+1}}=\Hb_{2^{d_1+d_2}} \otimes \Hb_2.$
  \end{enumerate}
\end{proposition}

Part (b) of Proposition~\ref{prop: hadamard} implies that $\bN={1 \over 2^{d_1+d_2}}\Hb\bS$, i.e., $n_{(\ba\bb)}={1 \over 2^{d_1+d_2}}\Hb_{(\overline{\ba\bb})}^T \bs$ where $\Hb_{(\overline{\ba\bb})}$ is the $(\overline{\ba\bb})$-th column of $\Hb$. With the above notation and transformation of variables, and by part (a) of Proposition~\ref{prop: hadamard}, the multinomial distribution in the contingency table \eqref{eq: likelihood_multinomial} can be written as
{\small
\begin{equation}\label{eq: likelihood_S}
 p(\bN=\bn|\bp)
 = {n! \over \prod_{\ba,\bb} n_{(\ba\bb)}!} \prod_{\ba,\bb}\bigg(\prod_{\ba',\bb'} p_{(\ba'\bb')}^{(-1)^{(\overline{\ba'\bb'})^T(\ba\bb)}}
 \bigg)^{s_{(\ba\bb)} \over 2^{d_1+d_2}}.
\end{equation}
}
We are now ready to introduce the interaction odds ratio (IOR):
\begin{definition}\label{def: ior}
We call $\lambda_{(\ba\bb)} =\prod_{\ba',\bb'} p_{(\ba'\bb')}^{(-1)^{\overline{(\ba'\bb')}^T(\ba\bb)}}$ the interaction odds ratio (IOR) with respect to the interaction $\dot{A}_{\ba}\dot{B}_{\bb}.$ Denote the vector of $\lambda_{(\ba\bb)}$'s by $\blambda$ and order the entries in the same way as in $\bE$.
\end{definition}
For each corresponding interaction, the IOR can be regarded as the ratio of the product of all white cell probabilities to the product of all blue cell probabilities. There are three cases for the IOR $\lambda_{(\ba\bb)}$:

(a) When $\ba=\zero$ and $\bb=\zero,$ $\lambda_{(\zero \zero)}= \prod_{\ba',\bb'} p_{(\ba'\bb')}.$ Note that the term $\lambda_{(\zero \zero)}^{n \over 2^{d_1+d_2}}$ does not involve $\bN$ and is constant.

(b) When $\ba=\zero$ but $\bb\neq \zero$ (or when $\bb=\zero$ but $\ba\neq \zero$),  then $\lambda_{(\ba\bb)}$ is a \emph{marginal interaction odds ratio} (MIOR) quantifying the balance in the marginal interaction variable $\dot{A}_{\ba}$ (or $\dot{B}_{\bb}$). For example, when $d_1=2$ and $d_2=1,$ $\lambda_{(110)}= { p_{(111)} p_{(110)} p_{(001)} p_{(000)}\over p_{(101)}p_{(100)}p_{(011)}p_{(010)}}$ which is related to the distribution of $\dot{A}_1\dot{A}_2.$ Note also that there are $2^{d_1} +2^{d_2}-2$ MIORs at depths $d_1$ and $d_2$.

(c) When $\ba \neq \zero$ and $\bb\neq \zero$, then $\lambda_{(\ba\bb)}$ is a \emph{cross interaction odds ratio} (CIOR) quantifying the balance in the cross interaction variable $\dot{A}_{\ba} \dot{B}_{\bb}.$ For example, when $d_1=2$ and $d_2=1,$ $\lambda_{(111)}={ p_{(111)} p_{(100)} p_{(010)} p_{(001)}\over p_{(110)}p_{(101)}p_{(011)}p_{(000)}}$ which is related to the distribution of $\dot{A}_1\dot{A}_2\dot{B}_1.$ Note also that there are $(2^{d_1}-1)(2^{d_2}-1)$ CIORs at depths $d_1$ and $d_2,$ which matches the degree of freedom for the $\chi^2$ test.

An important observation is that with the IOR, \eqref{eq: likelihood_S} becomes
{\footnotesize
\begin{equation}\label{eq: likelihood_IOR}
p(\bS=\bs|\blambda) = C_2(\bs) h_2(\blambda)\exp\bigg(\sum_{\ba \neq \zero} {s_{(\ba\zero)}\log\lambda_{(\ba\zero)} \over 2^{d_1+d_2}}+  \sum_{\bb \neq \zero} {s_{(\zero\bb)} \log\lambda_{(\zero\bb)} \over 2^{d_1+d_2}}+ \sum_{\substack{\ba \neq \zero \\ \bb \neq \zero}} {s_{(\ba \bb)} \log\lambda_{(\ba\bb)}\over 2^{d_1+d_2}}\bigg)
\end{equation}
}where $C_2(\bs)={n! \over \prod_{\ba,\bb} n_{(\ba\bb)}!} $ and $h_2(\blambda)=\lambda_{(\zero \zero)}^{n \over 2^{d_1+d_2}}.$ Therefore, we reparametrize the distribution in \eqref{eq: likelihood_multinomial} as a $(2^{d_1+d_2}-1)$-dimensional exponential family with log-IORs as natural parameters, and the symmetry statistics are \emph{complete sufficient statistics} for log-IORs. This fact is the basis of the binary expansion approach.

Similarly to the BID equations, we have a logarithm version of the BID equation:
\begin{theorem}\label{thm: be_log_ior}
Denote the vectors of the logarithm of entries in $\blambda$ and $\bp$ by $\blambda_l$ and $\bp_l$ respectively. We have $  \blambda_l= \Hb \bp_l$.
\end{theorem}

One important implication of \eqref{eq: likelihood_IOR} and Theorem~\ref{thm: be_log_ior} is that all information about dependence is contained in CIOR:
\begin{theorem}\label{thm: be_ior}
    $U_{d_1}$ and $V_{d_2}$ are independent if and only if $\lambda_{(\ba\bb)}=1$ for all CIORs.
\end{theorem}
Theorem~\ref{thm: be_ior} shows that the null hypothesis of the test \eqref{eq:test_be_ind} is equivalent to
{\small
\begin{equation}\label{eq: test_be_ior}
H_{0,d_1,d_2}: \text{ For all CIORs at depths }d_1 \text{ and } d_2, \lambda_{(\ba\bb)}=1.
\end{equation}
}
\vspace*{-\baselineskip}

\noindent We summarize the advantages of the reparametrization in \eqref{eq: likelihood_IOR} and the test \eqref{eq: test_be_ior}:\\

\vspace*{-0.7\baselineskip}

(a) Compared to the conventional parametrization in \eqref{eq: likelihood_multinomial}, the reparametrization in \eqref{eq: likelihood_IOR} is much more interpretable: Note that the cell probabilities in $\bp$ carry both marginal and joint information. On the other hand, the parametrization with $\blambda$ extracts all dependence information in CIORs and separates it from the marginal information in MIORs. Thus, CIORs are to contingency tables as correlations are to multivariate normal distributions. Tests of independence can therefore focus on CIORs, as we study in details in Section~\ref{sec: maxbet}.

(b) The sufficient statistics in the conventional parametrization are the cell counts $n_{(\ba\bb)}$'s, whose distribution is $Binomial(n,p_{(\ba\bb)}).$ This means that when $n$ is small, one often has $n_{(\ba\bb)}=0$ for many cells. These empty cells cause problems in the conventional tests. However, with the reparametrization \eqref{eq: likelihood_IOR}, the sufficient statistics $S_{(\ba\bb)}$'s instead have (after a linear transformation) a binomial distribution whose probability of success is the sum of $2^{d_1+d_2}-1$ cell probabilities. Therefore, by grouping the cells, $S_{(\ba\bb)}$'s provide much more information than $n_{(\ba\bb)}$'s and avoid the well-known problem of insufficient samples in many binning methods.

(c) Note that each cross interaction in the filtration corresponds to a unique CIOR, which measures some form of dependency. In Section~\ref{sec: maxbet}, we show that this consideration together with the number of CIORs $(2^{d_1}-1)(2^{d_2}-1)$ lead to an orthogonal decomposition of the $\chi^2$ test.

(d) The BID equation in Theorem~\ref{thm: bid} can be generalized for any three-way or multiway contingency table whose size is a power of 2. This fact allows extensions of the IOR reparametrization and the BET for testing independence of random vectors.

When the marginal distributions are unknown, for each observation $i$, we can similarly define $\hat{\dot{A}}_{k,i}=2\hat{A}_{k,i}-1$, $\hat{\dot{B}}_{k,i}=2\hat{B}_{k,i}-1,$ and $\hat{S}_{(\ba\bb)}=\sum_{i=1}^n \hat{\dot{A}}_{\ba,i}\hat{\dot{B}}_{\bb,i}$ for the cross interaction $\hat{\dot{A}}_{\ba}\hat{\dot{B}}_{\bb}.$ Now note the following simple corollaries from Theorem~\ref{thm: bid}: (a) $\bn_r$ and $\bn_c$ are invertible functions of $\hat{S}_{(\ba\zero)}$'s and $\hat{S}_{(\zero\bb)}$'s through a univariate BID equation, and (b) the bivariate sample BID equation holds for $\hat{\bS}$ and $\bn.$ With these facts, by using $\btheta$ and the proof of Theorem~\ref{thm: be_ior}, as well as conditioning on $\hat{S}_{(\ba\zero)}$ and $\hat{S}_{(\zero\bb)}$ in \eqref{eq: likelihood_multinomial}, we have
{\footnotesize
\begin{equation}\label{eq: likelihood_IOR_hypergeometric}
 p(\hat{\bS}_{(\ba\bb)}=\hat{\bs}_{(\ba\bb)} |\blambda_{(\ba\bb)},\hat{S}_{(\ba\zero)}, \hat{S}_{(\zero\bb)})
= C_2(\hat{\bs}_{(\ba\bb)}) h_3(\blambda_{(\ba\bb)}) \exp\bigg(\sum_{\substack{\ba \neq \zero \\ \bb \neq \zero}} {\hat{s}_{(\ba \bb)}\log\lambda_{(\ba\bb)} \over 2^{d_1+d_2}}\bigg)
\end{equation}
}
for some function $h_3(\blambda_{(\ba\bb)})$ as a normalizing constant.

Note that by conditioning on the counts of marginal interactions, the MIORs are eliminated, and we can focus on the CIORs for the analysis of dependence. Indeed, either by comparing \eqref{eq: cmvhypergeo} and \eqref{eq: likelihood_IOR_hypergeometric} or by the proof of Theorem~\ref{thm: be_ior}, we see that $\hat{U}_{d_1,i}$ and $\hat{V}_{d_2,i}$ are independent for each $i$ if and only if $\lambda_{(\ba\bb)}=1$ for all $\ba \neq \zero$ and $\bb\neq \zero$. Therefore, the tests of independence are unified in both of the cases of known and unknown marginal distributions to be \eqref{eq: test_be_ior}.

We remark here that reparametrization of the contingency table likelihood into odds ratios has been extensively studied in the past \cite{agresti1992survey}. The very recent paper \cite{ma2016fisher} also considered a factorization under the null hypothesis of independence. However, we are not aware of similar ideas of the connection to the Hadamard transform and the concept of IOR. Compared to existing analyses of contingency tables, the new reparametrization is more global to use all the observations. See a detailed discussion in Section~\ref{subsec: other}.

We also remark here that we are able to take advantage of the Hadamard transform only because the size of the contingency table is a power of $2$, which is a result of $\dot{A}_i$'s and $\dot{B}_j$'s in the binary expansions. If we were to take a different approach or to partition $[0,1]^2$ into different sizes, then we might not be able to have similar theory. This advantage is an important motivation of the binary expansion approach.

\section{The Max BET Procedure and Its Properties}\label{sec: maxbet}
\subsection{BET as an Multiple Testing Problem}\label{subsec: max}
In this section we return to the dependence detection problem, where we partition $[0,1]^2$ at the binary fractions based on Theorem~\ref{thm: be}. Therefore, the row and column total probabilities in the $2^{d_1}\times 2^{d_2}$ contingency table are $2^{-d_1}$ and $2^{-d_2}$ respectively when the marginal distributions are known, and the row and column total counts in the contingency table are $n2^{-d_1}$ and $n2^{-d_2}$ respectively when the marginal distributions are unknown and when $n$ is a multiple of $2^{\max\{d_1,d_2\}}.$

The discussions in Section~\ref{sec: bestat} suggest test statistics based on interactions $S_{(\ba\bb)}$ or $\hat{S}_{(\ba\bb)}.$ Direct application of the MLE of $\lambda_{(\ba\bb)}$ can result in similar disadvantages as $\chi^2$ tests as we discuss later. We instead construct a simple but optimal test statistic with the maximal symmetry statistics $\max|S_{(\ba\bb)}|$ or $\max|\hat{S}_{(\ba\bb)}|$ for $\ba \neq \zero$ and $\bb \neq \zero.$

The key observations of $S_{(\ba\bb)}$ are summarized below.
\begin{theorem}\label{thm: symmetry}
The following are equivalent:
\begin{enumerate}[(a)]
  \item $U_{d_1}$ and $V_{d_2}$ are independent.
  \item $\Eb[\dot{A}_{\ba}\dot{B}_{\bb}]=0$ for $\ba \neq \zero$ and $\bb\neq \zero$.
  \item $(S_{(\ba\bb)}+n)/2 \sim Binomial(n,1/2) $ for $\ba \neq \zero$ and $\bb\neq \zero$.
  \item $\Eb[S_{(\ba\bb)}]=0$ for $\ba \neq \zero$ and $\bb\neq \zero$.
  \item $\bE=\be_{\zero\zero}$ where $\be_{\zero\zero}$ is the $2^{d_1+d_2}$-dimensional standard basis $(1,0,\ldots,0)^T$.
\end{enumerate}
\end{theorem}

Note here that in Theorem~\ref{thm: symmetry}, the homogeneity in the distribution of $S_{(\ba\bb)}$ is due to the symmetry in $\dot{A}_{\ba}$ and $\dot{B}_{\bb}$ in the binary expansion. Indeed, the main intuition of Theorem~\ref{thm: symmetry} is the symmetry of independence: When $U_{d_1}$ and $V_{d_2}$ are independent, the counts of observations in the positive and negative regions should be similar for any cross interaction. On the other hand, when $U_{d_1}$ and $V_{d_2}$ are not independent, we expect some strong asymmetry between the numbers of points in white or blue.

When the marginal distributions are unknown, we have similar results on symmetry assuming $n$ is a multiple of $2^{\max\{d_1,d_2\}}$. When $\hat{U}_{d_1,i}$ and $\hat{V}_{d_2,i}$ are independent for each $i=1,\ldots,n$, the distribution of $(\hat{S}_{(\ba\bb)}+n)/4$ is $ Hypergeometric(n,n/2,n/2).$ An intuitive way to understand this is that if we assign all $n$ observations into a $2\times 2$ table according to $\hat{\dot{A}}_{\ba,i}=\pm 1$ and $\hat{\dot{B}}_{\bb,i}=\pm 1,$ $\hat{S}_{(\ba\bb)}$ is the difference in counts of the interaction $\hat{\dot{A}}_{\ba,i}\hat{\dot{B}}_{\bb,i}$ being $+1$ or $-1.$ We show below that the converse is also true.

\begin{theorem}\label{thm: symmetry_sample}
When $n$ is a multiple of $2^{\max\{d_1,d_2\}}$, the following are equivalent:
\begin{enumerate}[(a)]
  \item For each $i$, $\hat{U}_{d_1,i}$ and $\hat{V}_{d_2,i}$ are independent.
  \item $(\hat{S}_{(\ba\bb)}+n)/4 \sim Hypergeometric(n,n/2,n/2) $ for $\ba \neq \zero$ and $\bb\neq \zero$.
  \item $\Eb[\hat{S}_{(\ba\bb)}]=0$ for $\ba \neq \zero$ and $\bb\neq \zero$.
\end{enumerate}
\end{theorem}
Theorem~\ref{thm: symmetry} and Theorem~\ref{thm: symmetry_sample} reduce the test of independence to tests of marginal properties of $\Eb[S_{(\ba\bb)}]$ and $\Eb[\hat{S}_{(\ba\bb)}].$ In particular, these results show the equivalence between the BET at depths $d_1$ and $d_2$ and a multiple testing problem: The testing problems in \eqref{eq:test_be_ind} and \eqref{eq:test_be_ind_sample} are equivalent to testing if all cross interactions up to depths $d_1$ and $d_2$ are symmetric. The advantage of this consideration is two-folded: (a) We reduce the test of a joint distribution (difficult) to that of marginal ones (simple). (b) We reduce the test of dependence (difficult) to that of symmetry (simple).

Note that the equivalent multiple testing problem is about controlling the family-wise error rate (FWER): Rejecting any symmetry results in the rejection of independence. The simplest FWER control is the Bonferroni procedure, where the adjusted $p$-value is the minimum of 1 and the product of $(2^{d_1}-1)(2^{d_2}-1)$ and the smallest $p$-value of all marginal tests. We refer this procedure as the Max BET.

We illustrate the Max BET procedure at depths $d_1=2$ and $d_2=1$ with the 64 samples studied in Section~\ref{subsec: interaction}. The procedure consists of the following steps, as shown in Figure~\ref{fig: bid}:
\begin{enumerate}
\item[Step 1]: We count white and blue points for each cross interaction $\dot{A}_2\dot{B}_1$, $\dot{A}_1\dot{B}_1$, and $\dot{A}_1\dot{A}_2\dot{B}_1$ for $d_1=2$ and $d_2=1$.
\item[Step 2]: Among these three cross interactions, we look for the one with the strongest asymmetry, which is $\dot{A}_2 \dot{B}_1 $ with 25 in white and 39 in blue. The symmetry statistic is $S_{(011)}=-14$. The binomial $p$-value is $0.103$.
\item[Step 3]: Use the Bonferroni adjustment to multiply $3$ and get the overall $p$-value of the Max BET at depths $d_1=2$ and $d_2=1$ to be $0.310.$
\end{enumerate}

Would the Bonferroni procedure be overly conservative? Our observation is no because of the orthogonality of the symmetry statistics. A formal study of optimality of the Bonferroni procedure is in Section~\ref{subsec: optimality}. Here, we state some results on the joint properties of symmetry statistics which provide some intuition.
\begin{theorem}\label{thm: ortho}
\ \begin{enumerate}[(a)]
  \item When the marginal distributions are known and $U_{d_1}$ and $V_{d_2}$ are independent, the symmetry statistics $S_{(\ba\bb)}$'s are pairwise independent.
  \item When the marginal distributions are unknown and for each $i$, $\hat{U}_{d_1,i}$ and $\hat{V}_{d_2,i}$ are independent, $\hat{S}_{(\ba\bb)}$'s are uncorrelated.
  \item The classical $\chi^2$ test statistic $C$ is $C={1 \over n}\sum_{\ba\neq \zero,\bb\neq \zero} \hat{S}_{(\ba\bb)}^2.$
\end{enumerate}
\end{theorem}
Part (a) and (b) of Theorem~\ref{thm: ortho} imply that due to the orthogonality in the BID, each symmetry statistic provides non-redundant information. Furthermore, part (b) and (c) of Theorem~\ref{thm: ortho} imply that the $(2^{d_1}-1)(2^{d_2}-1)$ sample symmetry statistics $\hat{S}_{(\ba\bb)}$'s form an orthogonal decomposition of the $\chi^2$ test statistic whose degrees of freedom is also $(2^{d_1}-1)(2^{d_2}-1)$. Therefore, instead of aggregating the information through sum of squares in the $\chi^2$ statistic, we here take a divide-and-conquer approach. To follow up the discussions in Section~\ref{subsec: likelihood}, we summarize the advantages of our approach below and describe the details in Section~\ref{subsec: optimality} and Section~\ref{subsec: interpretation}.

(a) In \cite{arias2011global} and \cite{barnett2016generalized}, it was noted that when the number of hypotheses is large and the signals are rare and weak, using a Bonferroni type of multiple comparison control can substantially outperform $\chi^2$ tests. In our context, this means that when $d_1$ and $d_2$ are large and when the dependence is through only a few cross interactions, the $\chi^2$ test is ``wasting'' many degrees of freedom. Instead, using the Max BET can help discover weaker dependence.

(b) Interpretability. One major advantage of using cross interactions over the $\chi^2$ test is that the grouping arrangement of the white and blue cells for each interaction helps indicate the pattern of the dependence, as described earlier in Section~\ref{subsec: interaction}. When the dependence is through only a few of cross interactions, with the rejection of the Max BET, we can identify the strongest interactions between the variables. These strongest interactions can in turn help describe the dependence.

\subsection{Power and Optimality of the Max BET}\label{subsec: optimality}
In this section, we study the power of the Max BET when the marginal distributions are known. The uniform consistency of the Max BET at any depths $d_1$ and $d_2$ follows from classical analysis of contingency tables. Moreover, despite the conservative nature of the Bonferroni approach, we show below that the Max BET can be optimal in power for a large collection of alternative distributions:
\begin{theorem}\label{thm: optimality}
For any fixed $0<\delta<1/2,$ denote by $\cH_{1, d_1,d_2}^R$ the collection of alternative distributions $\Pb_{(U_{d_1},V_{d_2})}$ such that
{\small
\begin{enumerate}
\item $TV(\Pb_{(U_{d_1},V_{d_2})},\Pb_{0,d_1,d_2}) \ge \delta$;
  \item $  \|\bE-\be_{(\zero \zero)}\|_\infty \ge \sqrt{d_1+d_2}2^{-(d_1+d_2)/4}\|\bE-\be_{(\zero \zero)}\|_2.$
\end{enumerate}
}
Consider the testing problem
{\small
\begin{equation}\label{eq:test_be_inds}
H_{0,d_1,d_2}: \Pb_{(U_{d_1},V_{d_2})} = \Pb_{0,d_1,d_2} ~~v.s.~~ H_1: \Pb_{(U_{d_1},V_{d_2})} \in \cH_{1,d_1,d_2}^R.
\end{equation}
}
For large $d_1$ and $d_2$, we have the following:
\begin{enumerate}
  \item For any $\epsilon>0$, the Max BET with size $\alpha$ needs $n=O(2^{(d_1+d_2)/2}/\delta^2)$ samples to have power $1-\epsilon.$
  \item Let $\cT_\alpha$ be the collection of all measurable size-$\alpha$ tests: $\cT_\alpha=\{T_\alpha: \Pb_{0,d_1,d_2}(T_\alpha=1)\le \alpha\}$. If $n=o(2^{(d_1+d_2)/2}/\delta^2),$ then there $\exists 0<\epsilon'<1-\alpha$ such that
{\small      \begin{equation}
        \inf_{T_\alpha \in \cT_\alpha} \sup_{\Pb_{(U_{d_1},V_{d_2})} \in \cH_{1,d_1,d_2}^R} \Pb_{(U_{d_1},V_{d_2})}(T_\alpha=0)\ge 1-\alpha-\epsilon'.
      \end{equation}
}
\end{enumerate}
\end{theorem}

The magnitude of the minimal sample size requirement has been carefully studied in statistics, information theory and machine learning. It describes the minimal number of samples to uniformly detect certain departure from the independence and in turn indicates the uniform power of the test. Part 1 of Theorem~\ref{thm: optimality} states that such a requirement for Max BET is $O(2^{(d_1+d_2)/2}/\delta^2)$, which matches the optimal rate in \cite{paninski2008coincidence,acharya2015optimal}. Moreover, part 2 of Theorem~\ref{thm: optimality} asserts that if the sample size grows at any smaller rate, then for any test, there exist alternatives such that the power of this test is strictly bounded away from 1. In this sense, the Max BET is minimax in the sample size requirement.

Note that the consistency of $\chi^2$ tests is shown in \cite{agresti2011categorical,fienberg2007analysis} to require $n>2^{d_1+d_2}.$ This requirement is much higher than the magnitude $O(2^{(d_1+d_2)/2}/\delta^2)$ in Theorem~\ref{thm: optimality} and indicates that the power of $\chi^2$ test can be much less than that of the Max BET. One intuitive explanation of this fact is that $\chi^2$ tests rely on good estimates of each cell probability in the table, while in the Max BET $S_{(\ba \bb)}$'s are based on grouped cells to utilize all $n$ observations.

The condition $\|\bE-\be_{(\zero \zero)}\|_\infty \ge \sqrt{d_1+d_2}2^{-(d_1+d_2)/4}\|\bE-\be_{(\zero \zero)}\|_2$ compares the strongest signal to the overall signal in the space of alternatives and indicates the signals to take on a spiky form. It can also be regarded as (but is more general than) a sparsity constraint, as it can be satisfied when at most ${1 \over d_1+d_2}2^{(d_1+d_2)/2}$ (out of $(2^{d_1}-1)(2^{d_2}-1)$) cross interactions have non-zero means. Under this generalized form of sparsity, the Bonferroni approach is not overly conservative. In particular, Theorem~\ref{thm: optimality} is consistent with the results in \cite{arias2011global} under the ANOVA setting that when the signals are square-root sparse, the max test has better power than the $\chi^2$ test. Note also that such a condition over $\bE$ does not imply sparsity in $\bp.$ Therefore the optimal rate in \cite{paninski2008coincidence,acharya2015optimal} still applies and is attained by the Max BET.

The sample size requirement in Theorem~\ref{thm: optimality} also indicates that for a given sample size $n,$ one can expect to detect dependence up to a depth of about $\log_2 n.$ This result again explains the problem of non-uniform consistency: One cannot expect one test to uniformly detect all types of dependency, and with $n$ samples one can only reliably detect dependence up to a depth of about $\log_2 n$ in the binary expansion filtration approximation. Note again that with the $\chi^2$ test the depth can only go up to about ${1 \over 2} \log_2 n,$ which means it may not have good power for many forms of dependency.

\subsection{Interpretation of the Max BET}\label{subsec: interpretation}
In this section we explain the interpretations of the BET, i.e., we ask when the BET at depths $d_1$ and $d_2$ is rejected, where is the dependence? The BET can explain this question explicitly with the cross interactions, because it returns with the 50\% area with significantly more points..

\begin{figure}[hhhh]
\begin{center}
\includegraphics[width=\textwidth]{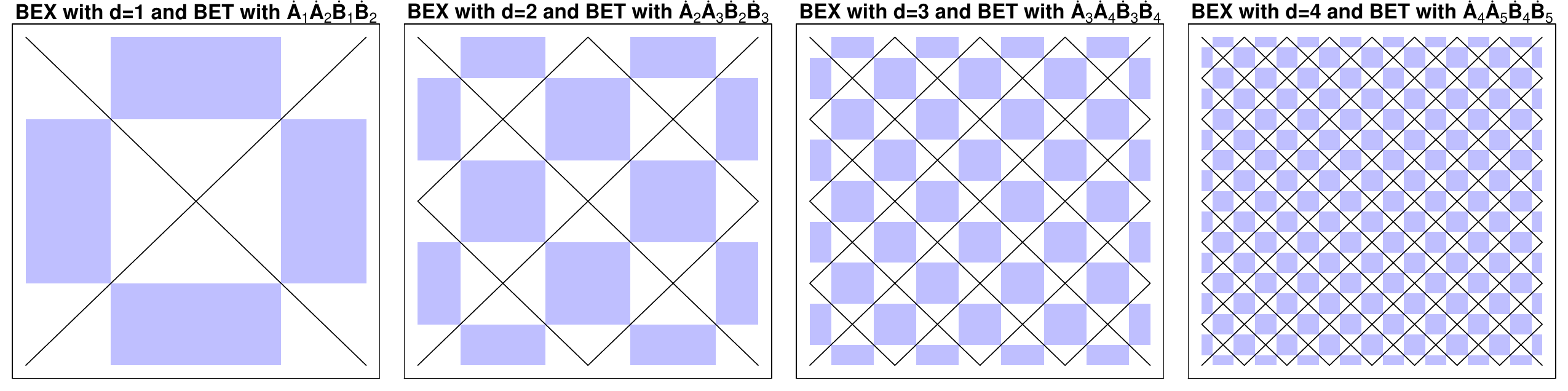}
\end{center}
\caption{
The bisection expanding cross (BEX) at $d=1,\ldots,4$ captured in the positive regions of the BET, which illustrates the interpretation of dependency in the BET.}\label{fig: betbex}
\end{figure}

We will explain some common patterns of dependence in simulation studies in Section~\ref{sec: numerical}. We will also illustrate the interpretation of BET with real data in Section~\ref{sec: stars} and Section~\ref{sec: tcga}. In what follows, we revisit the bisection expanding cross (BEX) as an example. See Figure~\ref{fig: betbex}. Note that with probability~1, samples of $(X_{d},Y_{d})$ on $BEX_d$ all fall in the positive region for $\dot{A}_d\dot{A}_{d+1}\dot{B}_d\dot{B}_{d+1}$. This is the strongest asymmetry of $BEX_d$, and the $p$-value for the Max BET at $d_1=d_2=d+1$ is $2 (2^{d+1}-1)^2/2^{n}$ which can be very small when $n$ is much larger than \red{$2d$}. Note that with the rejection of the Max BET at $d_1=d_2=d+1,$ the cross interaction $\dot{A}_d\dot{A}_{d+1}\dot{B}_d\dot{B}_{d+1}$ is also found to present the dependency between $X_d$ and $Y_d.$

With the above considerations, we explain the paradox following Proposition~\ref{prop: unif_bex}. For $(X_d,Y_d)$ on $BEX_d$, let $U_d$ and $V_d$ be the truncated variables in the marginal binary expansion of $X_d$ and $Y_d$ respectively. Note that $U_d$ and $V_d$ are independent. However, $U_{d+1}$ and $V_{d+1}$ are dependent, as is evidenced by the small $p$-value. These facts thus explain the seeming paradox: If we are at depths $d_1=d_2=d$, then the fact that $U_d$ and $V_d$ are independent implies that $X_d$ and $Y_d$ are $(d,d)$-independent, i.e., nearly independent. On the other hand, if we are at depths $d_1=d_2=d+1$, then the small $p$-value of the BET implies that $X_d$ and $Y_d$ are strongly non-independent. Therefore, being strongly non-independent or nearly independent depends on the choice of depth, and there is no contradiction in this example.

\subsection{Relations to Other Binning Methods}\label{subsec: other}
Although the binary expansion approach leads to multi-scale discretization, the BET is different from existing tests in the binning approach in several ways: (a) Many existing binning methods such as \cite{reshef2011detecting,kinney2014equitability} involve an optimization step in search of the optimal partition of data under some criteria such as mutual information. This step could be computationally expensive due to a search over many overlapping partitions which contain redundant information. Instead, the partitions based on interactions from the binary expansion filtration are created in a systematic manner with a natural hierarchy. The orthogonal design of interactions also saves much redundant information and improves the power. (b) Many binning tests may have problems of insufficient observations in small bins, while in the BET all $n$ samples are used repeatedly in an orthogonal manner which has advantages both for the level and power. (c) Many binning tests return a $p$-value based on permutations, which can again be computationally more expensive than the BET.

We also compare the Max BET with recent work in scan statistics \citep{walther2010optimal,ma2016fisher} which are based on rectangle scanning windows for local dependency. We note that some scanning method can be formulated in terms of the binary expansion statistics. For example, the FES in \cite{ma2016fisher} up to $(2,1)$-independence can be regarded as the following three tests of symmetry: $\Eb[\dot{A}_1\dot{B}_1]=0, \Eb[\dot{A}_2\dot{B}_1|\dot{A}_1=1]=0$ and $\Eb[\dot{A}_2\dot{B}_1|\dot{A}_1=-1]=0.$ Compared to the three tests of symmetry in the Max BET $\Eb[\dot{A}_1\dot{B}_1]=0$, $\Eb[\dot{A}_2\dot{B}_1]=0$ and $\Eb[\dot{A}_1\dot{A}_2\dot{B}_1]=0$, FES can be regarded as a conditional version of the BET. This conditional formulation can be advantageous in detecting local dependency, but may not have optimal power when the dependency is global and may have the insufficient sample problem discussed above. In the Max BET, the grouping of positive and negative regions does not necessarily result in a region of the rectangle shape but is more capable of detecting global dependency. Thus, each method has its advantageous scenarios.

\subsection{Issues in Practice}\label{subsec: sequential}
In this section we discuss issues of the Max BET that can happen in practice. The first issue is that we often do not know correct depths $d_1$ and $d_2$ where the dependency may be present. To address this issue, we propose a search over different depths and a second stage multiplicity control. This proposal is based on the observation that the approximation error in \eqref{eq: beuvd_joint} is $O_p(2^{- \min\{d_1,d_2\}}).$ Therefore, we can first test the hypotheses \eqref{eq: test_be_ior} for $d_1=d_2=d$ with $d=1,\ldots, d_{max}$, where $d_{max}$ reflects the desirable accuracy in the approximation. Then we can apply some further FWER multiplicity control procedure such as the Bonferroni method over the $d_{max}$ tests to ensure the overall FWER.

In practice, note that from \eqref{eq: beuvd_joint} $d_{max}=4$ provides good approximation to the true distribution. Note also that in order to avoid overlapping cross interactions in different depths, for each $d\ge 2$, one can test the symmetry of all added interactions involving $\dot{A}_d$ or $\dot{B}_d$, which are in $\sigma(U_{d},V_{d})$ but not in $\sigma(U_{d-1},V_{d-1}).$ We illustrate this procedure in Section~\ref{sec: numerical} and Section~\ref{sec: stars}. The effect of such multiplicity control on power is studied in Section 1.2 of the supplementary file.

Another practical issue for the empirical BET is that $n$ might not be a multiple of $2^{\max\{d_1,d_2\}}$, i.e., the column and row total counts might not be equal in the $2^{d_1} \times 2^{d_2}$ table. In this case, the reparametrization in Section~\ref{subsec: interaction} still applies, and the test for each cross interaction is still a Fisher'e exact test for $2 \times 2$ tables. However, the distribution of a symmetry statistic (after a linear transformation) is not necessarily $Hypergeometric(n,n/2,n/2).$ In general, instead of $n/2$'s, the parameters for the hypergeometric distribution are numbers of observations for which the marginal interactions are positive. Thus, symmetry and homogeneity might be lost in this case. Nonetheless, the BET still applies for any sample size $n \ge 2^{\max\{d_1,d_2\}}$ (otherwise there exist cross interactions for which all observations are positive). Moreover, when $n$ is large, one can use the normal approximation in \cite{kou1996asymptotics} for these tests.

\section{Connection to Computing}\label{sec: computing}
The binary expansion approach is partially motivated by its close connections to the current computing system, which is based on a binary architecture. By turning an electrical circuit ``on'' (represented by ``1'') and ``off'' (represented by ``0''), computers process information with unprecedented speed and power. In particular, each decimal number in computing is processed as a rounded version of its binary representation. For example, calculations of $0.1_{10}=0.000110011\ldots_2$ are based on a rounded version of $0.000110011\ldots_2$ to certain bits (depending on a 32-bit or 64-bit computing system).

The key observation here is that {\it the binary representation of a decimal number is precisely its binary expansion!} The $\{A_k\}_{k=1}^{d_1}$ and $\{B_k\}_{k=1}^{d_2}$ in the BEStat approach directly correspond to the first $d_1$ and $d_2$ bits of $U$ and $V$ respectively in current computing systems. This fact implies that as long as a statistician is processing data with a computing device (desktop, laptop, smartphone, hand-held calculator...), the $\{A_k\}_{k=1}^{d_1}$ and $\{B_k\}_{k=1}^{d_2}$ are given to him/her automatically. These binary bits are hidden resources of data available for statisticians from computers. We often use bits for computing, but {\it bits are data!} We can construct statistics and make inference with bits, and the BET at depths $d_1$ and $d_2$ can be explicitly interpreted as testing whether the data are independent up to the first $d_1$ and $d_2$ bits.

Moreover, the BEStat approach provides statisticians the access to the most fundamental level of the computing system and enables direct operations over bits. For example, the cell locating process of a data point in the contingency table can be done through some bitwise Boolean operations over the $a_k$'s and $b_k$'s. Such bitwise operations are known to be computationally efficient. We develop such a bitwise algorithm of the BET in a separate paper \citep{zhao2019fast}, where the procedure is shown to improve the speed of existing methods by orders of magnitude.

\section{Simulation Studies}\label{sec: numerical}
In this section, we use simulation studies to compare the Max BET and existing nonparametric methods. For the Max BET, we consider the empirical CDF transformation and consider the second stage multiplicity control over depths with the Bonferroni procedure with $d_{max}=4$, as discussed in Section~\ref{subsec: sequential}. For comparison, we consider the Hoeffding's D test from the CDF approach, the distance correlation from the distance approach, the default KNN-MI method from the binning approach, and the very recent method of FES. We consider the $\chi^2$ test for the same contingency table for the Max BET with $d_1=d_2=4$ too.

We compare the power the above methods over common dependency structures such as linear, parabolic, circular, sine, and checkerboard, which are widely considered in evaluation of tests of independence \citep{reshef2011detecting,heller2012consistent,kinney2014equitability,filippi2015bayesian}. We also consider the local dependency setting in \cite{ma2016fisher}. The scenarios are designed by adapting those in \cite{ma2016fisher} with an emphasis on small sample performance with a fixed sample size $128$. The level of the tests are set to be $0.1$. We simulate each of the scenarios at 10 different noise levels to present the whole range of power. The details of the setting are summarized in Table~\ref{table: power}.
\begin{table}[hhhh]
\def\arraystretch{1}%
\centering
{\footnotesize
\begin{tabular}{l|l|l}
Scenario & Generation of $X$ & Generation of $Y$\\
\hline
Linear & $X=U~~~~~$ & $ Y=X+6\epsilon$ \\
Parabolic & $X=U~~~~~$ & $  Y=(X-0.5)^2+1.5\epsilon$\\
Circular & $X=\cos \vartheta+2.5\epsilon~~~~~$ & $  Y=\sin\vartheta+2.5\epsilon'$\\
Sine & $X=U~~~~~$ & $  Y=\sin(4 \pi X)+8 \epsilon$\\
Checkerboard & $X=W+\epsilon~~~~~$ & $  Y=\left\{\begin{array}{ll}
V_1 +4\epsilon' & \text{ if }  W=2\\
V_2 +4\epsilon'' & \text{ otherwise}
\end{array}\right.$\\
Local & $X=G_1~~~~~$ & $  Y=\left\{\begin{array}{ll}
X+\epsilon &  \text{ if } 0 \le G_1 \le 1 \text{ and } 0 \le G_2 \le 1\\
G_2 & \text{ otherwise}
\end{array}\right.$\\
\hline
\end{tabular}
}
\caption{Simulation scenarios: At each noise level $l=1,\ldots,10$, $\epsilon, \epsilon', \epsilon'' \stackrel{iid}{\sim} \cN(0,(l/40)^2)$, and the following variables are all independent: $U \sim Uniform[0,1]$, $\vartheta \sim Uniform[-\pi,\pi]$, $W \sim Multi-Bern(\{1,2,3\}, (1/3,1/3,1/3))$, $V_1 \sim Bern(\{2,4\},(1/2,1/2))$, $V_2 \sim Multi-Bern(\{1,3,5\}, (1/3,1/3,1/3))$, $G_1, G_2 \stackrel{iid}{\sim} \cN(0,1/4)$.}\label{table: power}
\end{table}

\begin{figure}[hhhh]
\begin{center}
\includegraphics[width=\textwidth]{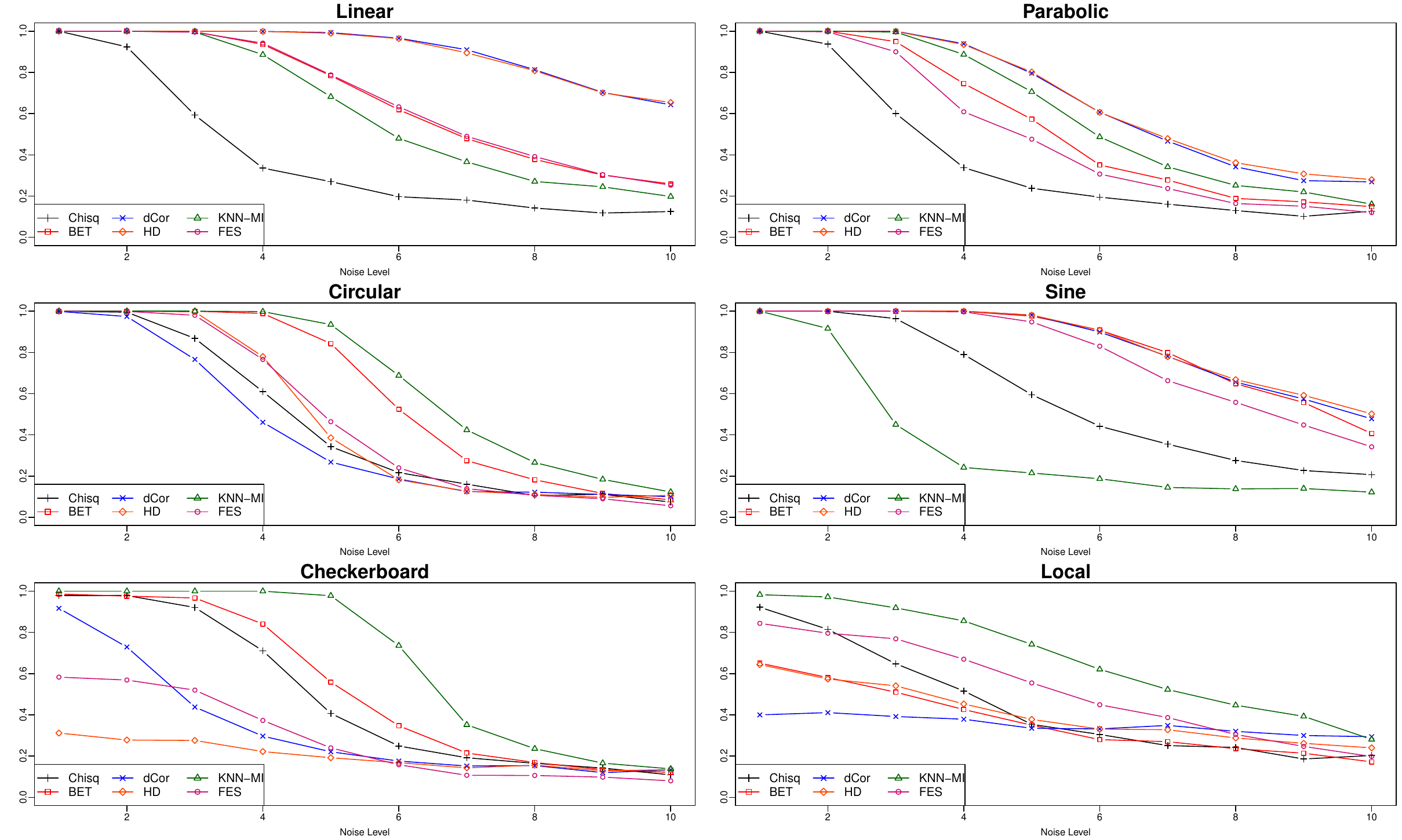}
\end{center}
\caption{
Comparison of powers from six nonparametric tests of independence: the two-stage Max BET with empirical CDF and with $d_{max}=4$ (BET), $\chi^2$ test for the discretization when $d_1=d_2=4$ (Chisq), distance correlation (dCor), Hoeffding's D (HD), $k$-nearest neighbor mutual information (KNN-MI), and Fisher exact scanning (FES).}\label{fig: bet_simu_power}
\end{figure}
The power curves of the six nonparametric tests of independence are presented in Figure~\ref{fig: bet_simu_power}. Generally speaking, as is found similarly in \cite{ma2016fisher} and many other papers, no test can uniformly dominate all others in all settings. In what follows, we separate the detailed discussions of the first five scenarios (linear, parabolic, circular, sine, and checkerboard) and the last scenario (local).

\begin{figure}[hhhh]
\begin{center}
\includegraphics[width=\textwidth]{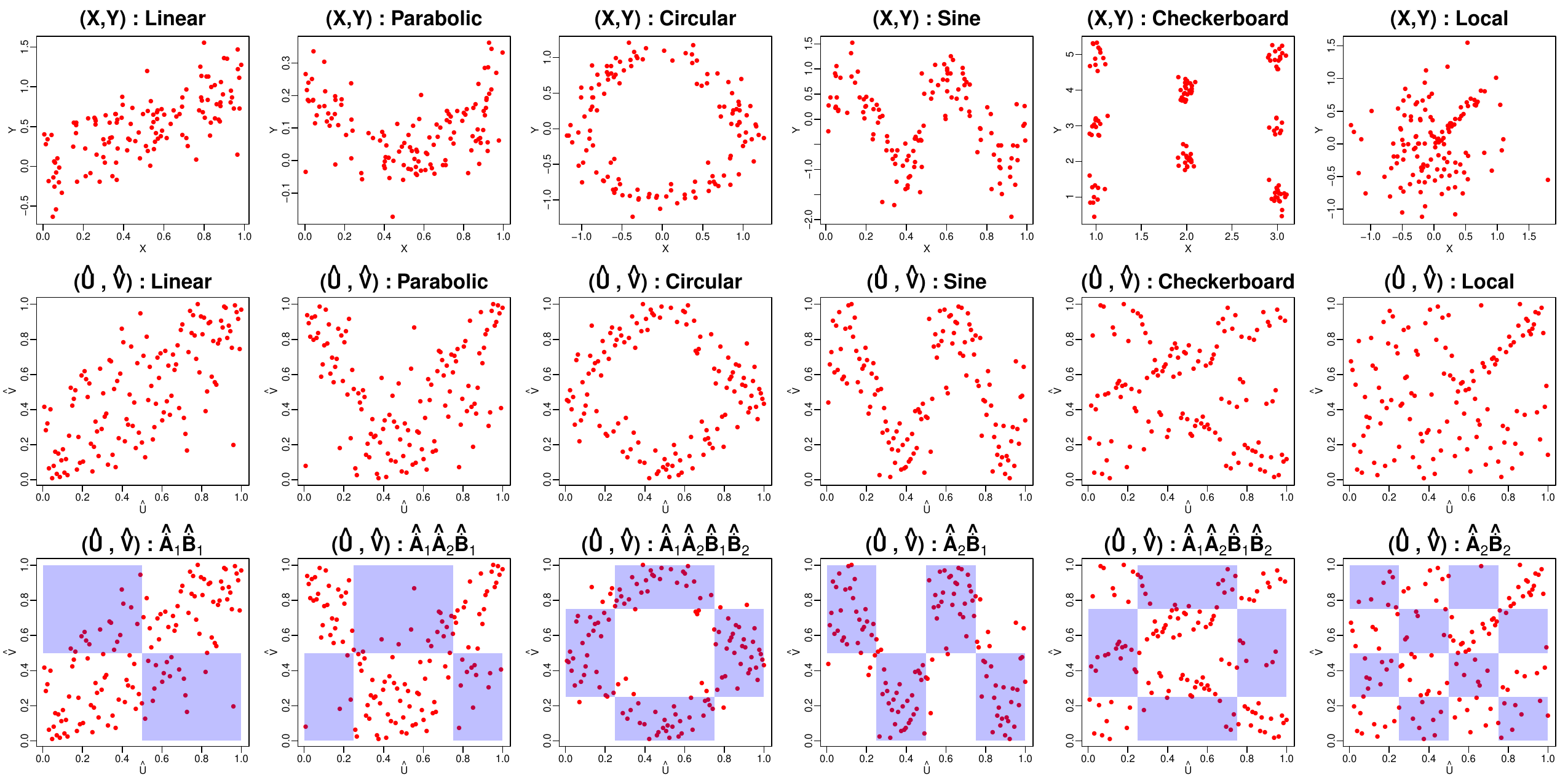}
\end{center}
\caption{
The BET interpretations of dependency patterns. The observations are generated as in Table~\ref{table: power} with noise level $l=2$. The first row shows the scatterplots of original data $(X,Y)$. The second row shows the corresponding empirical copula $(\hat{U},\hat{V})$ for $i=1,\ldots,128.$ The third row shows the cross interaction of the strongest asymmetry, which the BET returns with the rejection of independence null.}\label{fig: bet_simu_interpretation}
\end{figure}

In the first five scenarios where the dependency is global, we notice that each existing method has shown some limitations: In the linear and parabolic setting, the $\chi^2$ test provides the least power. In the circular setting, distance correlation provides the least power. In the sine setting, KNN-MI provides the least power. In the checkerboard setting, Hoeffding's D and FES provide the least power, which is partially due to the fact that observations in this setting are locally independent. On the other hand, the BET never provides the least power under these common relationships. One reason of such robustness of the BET is that the global dependency in these settings can be well explained through only a few cross interactions in the binary expansion, as can be seen in Figure~\ref{fig: bet_simu_interpretation} and in discussions below. Therefore, the minimaxity in Theorem~\ref{thm: optimality} guarantees that the BET has reliable power against a large class of alternative distributions. We also note here that to echo with the discussions in Section~\ref{subsec: other}, the BET has better power than FES in most of these global dependency settings because of its global grouping of cells. On the other hand, FES has better performance in the local dependency setting, as we discuss below.

We now turn to the setting of the local relationship. The BET does not perform well because observations in this setting are independent outside the area with the local dependency. Therefore, the global grouping of cells in the BET does not provide more information than a few local cells. In this case, the condition in Theorem~\ref{thm: optimality} can be violated as many cross interactions are asymmetric with weak signals. As shown in Figure~\ref{fig: bet_simu_power}, this limitation of the BET can be remedied by scanning based binning methods such as FES, which focuses on local dependency, or clustering based binning methods such as KNN-MI, which performs well on mixtures of distributions.

One useful property of the BET is its interpretability of dependency based on the interactions of binary variables, which we illustrate in Figure~\ref{fig: bet_simu_interpretation}. In each column, we present a simulated dataset in each scenario with noise level $l=2$. In the first five scenarios, the global dependency in the data is well explained by a corresponding cross interaction: Observations with linear dependency tend to fall in the positive region of $\hat{\dot{A}}_1\hat{\dot{B}}_1$, observations with the parabolic dependency tend to fall in the positive region of $\hat{\dot{A}}_1\hat{\dot{A}}_2\hat{\dot{B}}_1$, observations with circular dependency tend to fall in the negative region of $\hat{\dot{A}}_1\hat{\dot{A}}_2\hat{\dot{B}}_1\hat{\dot{B}}_2$, observations with the sine dependency tend to fall in the negative region of $\hat{\dot{A}}_2\hat{\dot{B}}_1$, observations with the checkerboard dependency tend to fall in the positive region of $\hat{\dot{A}}_1\hat{\dot{A}}_2\hat{\dot{B}}_1\hat{\dot{B}}_2.$ Since these common global dependency patterns can be well explained by a single cross interaction, Theorem~\ref{thm: optimality} applies and the Max BET has good performance in terms of power as shown in Figure~\ref{fig: bet_simu_power}.

The local dependency in the last scenario is also well captured by the positive region of $\hat{\dot{A}}_2\hat{\dot{B}}_2$, particularly in the four upper right cells. However, outside this region the variables are independent, so the interpretation of dependency is rather explained by a local and conditional cross interaction $\hat{\dot{A}}_2\hat{\dot{B}}_2$ given $\{\hat{\dot{A}}_1=1,\hat{\dot{B}}_1=1\}$, than by the global cross interaction $\hat{\dot{A}}_2\hat{\dot{B}}_2$. In this case, scanning based binning methods such as FES provide better interpretation of the local dependency.

\section{Are Stars Randomly Distributed in the Sky?}\label{sec: stars}
In this section we study the curious question of whether stars in the night sky are randomly distributed. Despite a simple statement of this long standing question, we are not aware of any complete scientific theory that explains the phenomenon with a confirming or disconfirming answer. In what follows, we provide some statistical analysis of this problem.

To study this question, we collected the galactic coordinates of the 256 brightest stars in the night sky \citep{perryman1997hipparcos}. The galactic coordinates are essentially spherical coordinates with the Sun as the center. These coordinates consist of radius, longitude $\phi \in [0,2\pi)$ and latitude $\varphi \in (-\pi/2,\pi/2]$. We ignore the radius information and focus on the unit sphere. Since the density of the uniform distribution over the unit sphere is proportional to $\cos{\varphi} d\phi d \varphi$, as long as $X=\phi$ and $Y=\sin{\varphi}$ of the stars are independent, the stars are uniformly distributed in the night sky.

\begin{figure}[hhhh]
\begin{center}
\includegraphics[width=\textwidth]{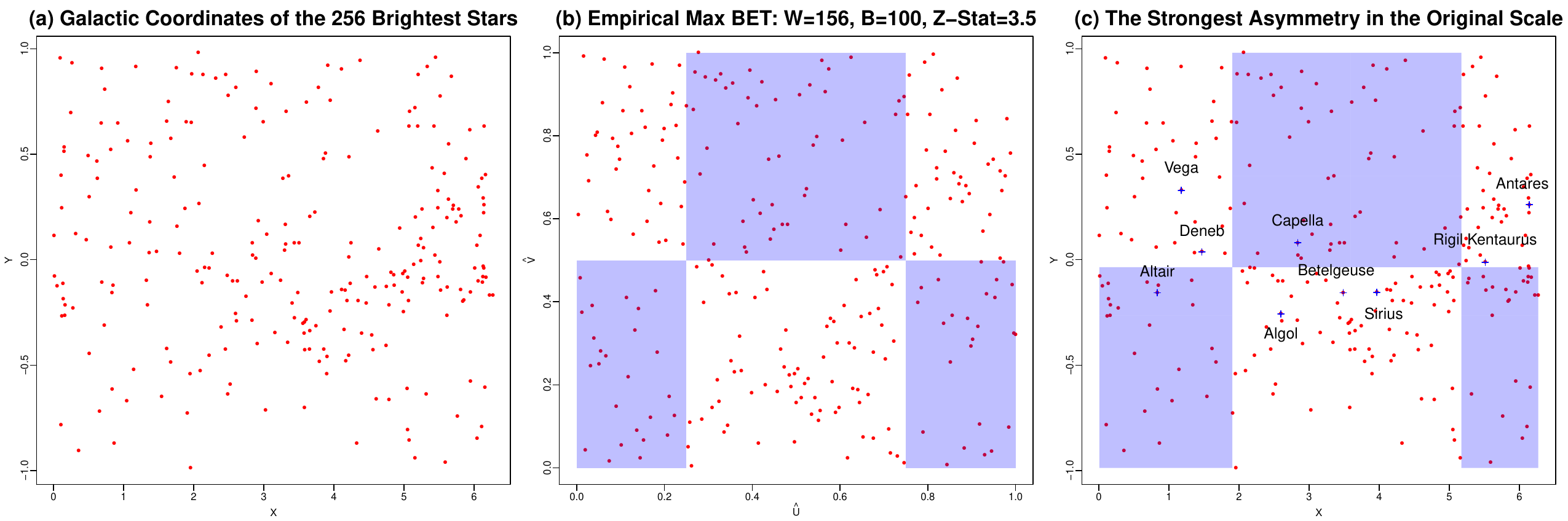}
\end{center}
\caption{(a) The longitude and sine latitude of the 256 brightest stars in the night sky. (b) The strongest asymmetry for the BET at $d=2$ is found to be the interaction $\hat{\dot{A}}_1\hat{\dot{B}}_1\hat{\dot{A}}_2$. (c) The strongest asymmetry in the original scale and some famous stars along the Milky Way.}\label{fig: stars}
\end{figure}
We first consider some classical tests of independence. The sample correlation between $X$ and $Y$ is $-0.07$ with a $p$-value of $0.264$, which is not significant. The distance correlation between $X$ and $Y$ is $0.137$ with a $p$-value of $0.064.$  Hoeffding's D test returns with a $p$-value of $0.103$. These $p$-values indicate some evidence against independence. The KNN-MI test provides a $p$-value of $0.02$, which is strong evidence against independence. However, this $p$-value does not provide any information about the relationship between $X$ and $Y$, and the dependence pattern is still unclear even when we rejected the null.

We now consider applying the two-stage empirical Max BET with $d_{max}=4$ on these data. The BET returns the strongest asymmetry $\hat{\dot{A}}_1\hat{\dot{A}}_2\hat{\dot{B}}_1$, where 156 stars are in the positive region and 100 are in the negative region. Thus, $\hat{S}_{(111)}=56$ and the approximate $z$-statistics is $3.5$ with the overall $p$-value $0.019$. Besides the strong evidence against independence, one important advantage of the BET is that we can also visualize the dependency upon rejection. In part (c) of Figure~\ref{fig: stars}, we transform the interaction in part (b) back to the original scale. Note that the labeled stars are well-known to be along the Milky Way in the night sky. Indeed, the Milky Way in the night sky is where stars in the galaxy cluster together, and its shape is captured by the positive region of $\dot{A}_1\dot{A}_2\dot{B}_1$. This fact explains the dependency in this data and the significance of the BET.

We note here that the application of FES to the star data returns with a $p$-value of $0.032$ with the strongest local dependency in $\hat{\dot{A}}_2\hat{\dot{B}}_1$ given $\{\hat{\dot{A}}_1=1\}.$ Compared with the BET which uses all $256$ observations to detect the dependency in $\hat{\dot{A}}_1\hat{\dot{A}}_2\hat{\dot{B}}_1$, the $p$-value of FES is higher because it only uses $128$ observations in the detection of local dependency when $\{\hat{\dot{A}}_1=1\}.$ In terms of interpretation, the FES only explains the dependency in the data with the ``right arm'' of the milky way, whereas the BET captures the entire milky way with an global cross interaction $\hat{\dot{A}}_1\hat{\dot{A}}_2\hat{\dot{B}}_1$.

A caveat here is that we regard the above analysis more as an illustration of the BET method rather than a scientific discovery, which requires a much more careful study. For example, the only strong assumption in the BET approach is the i.i.d. assumption on the observations. This assumption might be violated when the data points are stars. Moreover, we also note that the radius, which is excluded from this study, plays an important role in the location of stars. However, the interpretations from the BET can still be of immediate practical value: For example, it can help people find bright stars in the night sky.

\section{Exploratory Data Analysis of TCGA Data}\label{sec: tcga}
\subsection{Nonlinearity and Mixture of Subtype Distributions}
Conventional exploratory data analysis (EDA) of small multivariate datasets usually starts with a scatterplot matrix, see \cite{buja1992computing} and \cite{cleveland1993visualizing} for good reviews. Pairwise scatterplots can help people find interesting dependency patterns among variables, which can in turn suggest further statistical or scientific investigation. However, for high-dimensional data, the scatterplot matrix is not feasible since there are too many pairwise plots to inspect \citep{sun2014putting}. Common EDA tools in this situation such as principal component analysis, can only show high-level structure in the data, and focus mainly on linear relationships of variables. The BET can provide an alternative approach for such EDA due to the interpretability of its $p$-value. We illustrate this idea below in the context of breast cancer classification.

The TCGA lobular freeze breast cancer data in \cite{cancer2012comprehensive} and \cite{ciriello2015comprehensive} contain gene expression intensities of 817 subjects, about $2/3$ of which, or 544 samples, are used here as a training set and the remaining 273 observations are used as a test set. This dataset is based on 16615 genes. There are five subtypes groups indicated in this dataset. In what follows, we focus on basal-like breast cancer, which is known to be more aggressive, more difficult to treat, and have poorer prognosis compared to the other subtypes \citep{perou2000molecular}. Accurate classification of this subtype of breast cancer is thus very important for the health quality of patients.

The goal of this analysis is to use the BET as an EDA tool in the training dataset in search for nonlinear dependency between pairs of genes. Once a pair is identified in the EDA phase we look in the literature for mentions of the two corresponding genes and study their connection to subgroup typing. We also use the test dataset for confirmatory analysis.

Why can nonlinear dependency be related to studies on subgroup typing? As we illustrate below, one source of nonlinearity could be mixture of different subtype distributions. Intuitively, some genes might have different joint behavior under different subtypes of cancer. Such distributional differences could be in location, scale, covariance and other moments. When these different bivariate distributions are mixed together, some nonlinear dependency pattern could be created in the pooled joint distribution. Since the BET can capture nonlinear dependency patterns and indicate the form of nonlinearity, once a pair is identified by the BET, we hope to track back with the label information to find interesting pairs of genes that are related to different subtypes of breast cancer.

We first prepare the data by excluding genes which had non-unique entries in intensities. Such ties are results of the thresholding step in the preprocessing, and we exclude these genes here for simplicity. This filtering step results in 10107 genes in the remaining data. In the EDA phase with the training dataset, we scan over all pairs of these 10107 genes with the BET based on the empirical CDF transformation and depths $d_1=d_2=2$, and the $p$-value are calculated based on the large sample normal approximation of hypergeometric distribution in \cite{kou1996asymptotics}. This approach leads to a total of ${10107 \choose 2} =51070671 \approx 5 \times 10^7$ comparisons. We control the multiplicity over these comparisons through the Bonferroni method. We use the level $0.1$ threshold for multiplicity adjusted $p$-values to determine whether a pairing is interesting enough to follow up in the literature.

We emphasize here that many existing nonparametric dependence detection methods, such as Hoeffding's D, distance correlation, KNN-MI and FES, are not suitable for this EDA task for the following reasons:

(a) Classical methods such as Hoeffding's D, distance correlation, and KNN-MI do not provide clear interpretation upon rejection of independence. For example, even if the tests based on them are significant, they cannot distinguish pairs of genes with nonlinear dependency from pairs of genes with linear dependency.

(b) Although mutual information based methods such as KNN-MI have good power against mixtures of distributions, the $p$-value of KNN-MI is obtained through permutations. With the Bonferroni control over $5 \times 10^7$ pairwise tests, we need at least $5 \times 10^8$ random permutations for each test in order to have a valid significance level of $0.1$. The computational expense is prohibitive.

(c) Although FES provides interpretation of local dependency, it does not allow users to specify a global form of dependency in search of interesting relationships between variables. Thus it cannot identify pairs of genes with global nonlinear dependency. See the discussions below.

\subsection{Results from TCGA Data}
In the EDA phase, the BET rejects independence over more than $10000$ pairings of genes out of $5 \times 10^7$. Out of these pairs of genes, we can focus on some particular form of dependency. For example, we can restrict on pairs of genes whose dependency can be explained by the cross interaction $\hat{\dot{A}}_1\hat{\dot{A}}_2\hat{\dot{B}}_1\hat{\dot{B}}_2$. This consideration results in only 84 pairs of genes. Note that this specification process of global dependency is not possible with FES, nor other existing methods.
\begin{figure}[hhhh]
\begin{center}
\includegraphics[width=\textwidth]{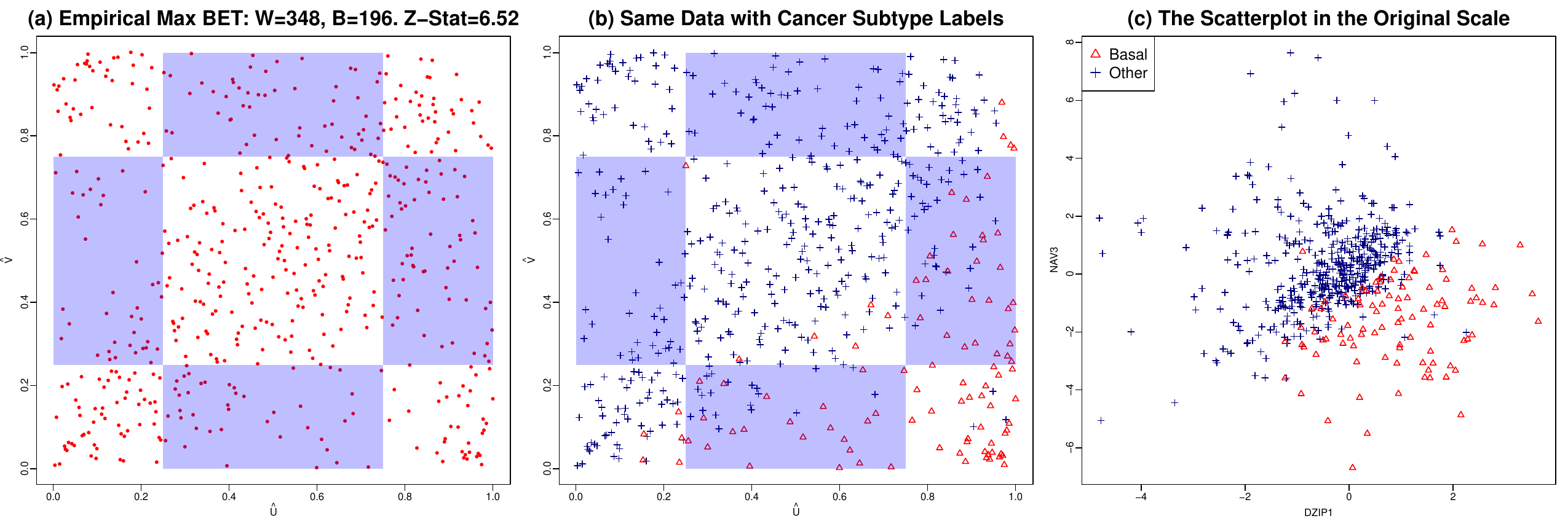}
\caption{(a) The BET with $d=2$ for two genes in the TCGA data. There are 348 observations and 196 observations in the empirical copula distribution falling into the positive and negative regions of $\hat{\dot{A}}_1\hat{\dot{A}}_2\hat{\dot{B}}_1\hat{\dot{B}}_2$ respectively. The $z$-statistic of the difference is 6.52. (b) The same two genes with the labels shown. Basal-like breast cancer patients are marked with a red triangle. (c) The scatterplot of same gene expressions in the original scale.}\label{fig: tcga}
\end{center}
\end{figure}
Of those 84 pairs of $\hat{\dot{A}}_1\hat{\dot{A}}_2\hat{\dot{B}}_1\hat{\dot{B}}_2$ dependency, we focus on \verb|DZIP1| and \verb|NAV3|. For this pair of genes, there are 348 observations and 196 observations falling into the positive and negative regions of $\hat{\dot{A}}_1\hat{\dot{A}}_2\hat{\dot{B}}_1\hat{\dot{B}}_2$ respectively. See Figure~\ref{fig: tcga}(a). The symmetry statistic is $\hat{S}_{(1111)}=152$, and the $z$-statistic of the difference is $6.52$, making the $p$-value of the BET to be $6.5\times 10^{-10}$. After multiplying $5\times 10^7$ for the Bonferroni control, the overall adjusted $p$-value is $0.033,$ which is strong evidence against the independence null. Furthermore, from the interaction $\hat{\dot{A}}_1\hat{\dot{A}}_2\hat{\dot{B}}_1\hat{\dot{B}}_2$ we could see interesting dependency patterns: In part of the data there exists strong monotone increasing dependency, while there is a cluster of observations above the third quartile of $U$ and below the first quartile of $V$. These patterns make the overall dependency nonlinear, which is captured by $\hat{\dot{A}}_1\hat{\dot{A}}_2\hat{\dot{B}}_1\hat{\dot{B}}_2$.

The above EDA with the BET suggests an interesting question: What is the reason of this nonlinear dependency? By adding the label of basal-like breast cancer, the cluster of observations in the lower right white box can be explained as a result of the joint distribution of the two genes under this subtype. From Figure~\ref{fig: tcga}(b) we see clearly that basal-like breast cancer patients tend to have higher \verb|DZIP1| intensity and lower \verb|NAV3| intensity. We also make the scatterplot of the same two genes in the original scale in Figure~\ref{fig: tcga}(c), and we see that the bivariate distribution of \verb|DZIP1| and \verb|NAV3| under the basal-like subtype has different location and scale and is almost disjoint from the rest of the data. This fact explains the reason of nonlinearity in the pooled distribution: When the bivariate distribution of this subtype is mixed together with those of other subtypes, some nonlinearity pattern is created. With the identification of this nonlinearity from the BET and with the label information, we can retrospectively extract such mixtures of different subtype distributions.

By searching the medical literature, we find both genes have been individually investigated and are confirmed to be highly related to basal-like breast cancer. For examples, the relationship between \verb|DZIP1| and basal-like cancer is studied by \cite{Kikuyama2012204,ShigunovShigunov2014}, and similar studies for \verb|NAV3| are done in \cite{maliniemi2011,cohen2015navigator}. However, we are not able to find results on the joint behavior of these two genes. The BET result indicates that this joint behavior could be scientifically important, as these two genes behave dramatically different under the basal-like subtype. This further suggests the possible existence of some biological functional relationships between these two genes and this subtype of cancer. This could be an interesting issue to investigate.

\subsection{Improvements in Classification}
Statistically, the above EDA with the BET suggests that \verb|DZIP1| and \verb|NAV3| could jointly be good predictors of basal-like breast cancer. We validate this conjecture with the test dataset of 273 subjects. We use the $k$-nearest neighbor classification method with $k=1$. The classification accuracy in the test dataset is 91\%. We assess this performance with cross-validation and observe similar results. Note that if we were to use \verb|DZIP1| or \verb|NAV3| alone for the classification task, the accuracy was 79\% and 76\% respectively, i.e., each of them is a good predictor but far from perfect. However, by combining these two genes and using the joint distribution for classification, we substantially improve the classification accuracy.

Existing classification studies are usually based on a selected set of many variables. One drawback of such studies is lack of interpretability. With some black box selection procedure over many variables, the effect of each variable is hard to scientifically interpret. On the other hand, the BET analysis can help identify pairs of variables which have high potential joint classification power, and explanations of the effects of variables can be obtained from the pattern of the nonlinear dependency. Therefore, the BET can be a useful EDA tool in practice: It provides $p$-values that we can see.

\section{Summary and Discussions}\label{sec: summary}
Nonparametric dependence detection is an important problem in statistics. To avoid the power loss due to non-uniform consistency, we introduce the concept of binary expansion statistics (BEStat), which combines four classical statistical wisdoms: copula, filtration, orthogonal design and multiple testing. The proposed binary expansion testing (BET) framework combines the strength from these wisdoms and enjoys the invariance property from the copula distribution, universality, identifiability and uniformity from the filtration, orthogonality and symmetry from the orthogonal design, and interpretability from multiple testing. The binary expansion approach also facilitates efficient bitwise computing implementation.

Two important potential generalizations are nonparametric tests of independence for general categorical variables and for random vectors. For general contingency tables, the filtration and the separation of marginal and joint information need to be developed carefully. For random vectors, the binary expansion filtration approximation in \eqref{eq: beuvd_joint}, the BID equation in Theorem~\ref{thm: bid} and the IOR reparametrization can all be generalized. We welcome further thoughts on related topics for deeper understanding of dependence and useful procedures in practice.

\section*{Supplementary Materials}
Online supplementary materials for this article include additional numerical studies, proofs of the results, and R functions used in the numerical studies.

\section*{Acknowledgements}
The author thanks Richard Berk, Larry Brown, Andreas Buja, Edward George, Arun Kumar Kuchibhotla, Linda Zhao and Zhigen Zhao for inspiring discussions that stimulated this research. The author also thanks Edoardo Airoldi, Mike Baiocchi, Shankar Bhamidi, Bhaswar Bhattacharya, Rong Chen, Jessi Cisewski, Bradley Efron, Jianqing Fan, Dean Foster, Andrew Gelman, Jan Hannig, Ruth Heller, Peter Hoff, Katherine Hoadley, Xiaoming Huo, Pierre Jacob, Vinay Kashyap, Michael Kosorok, S.C. Samuel Kou, Duyeol Lee, Michael Levine, Ping Li, Yun Li, Xihong Lin, Oliver Linton, Han Liu, Jun Liu, Mike Love, Li Ma, Zongming Ma, Steve Marron, Xiao-Li Meng, Joel Parker, Charles Perou, Vladas Pipiras, Richard Samworth,  David Siegmund, Dylan Small, Robert Stine, William Strawderman, Weijie Su, G\'{a}bor Sz\'{e}kely, Xinlu Tan, Yihong Wu, Han Xiao, Daniel Yekutieli, Ming Yuan, Yuan Yuan, Cun-Hui Zhang, Nancy Zhang, Tingting Zhang, and Harrison Zhou for valuable comments and suggestions. This research is completed while the author is visiting Princeton University. The author thanks Jianqing Fan and the Department of Operations Research and Financial Engineering at Princeton for the warm hospitality.

\section*{Funding}
This research is partially supported by NSF DMS-1613112 and NSF IIS-1633212.

\bibliographystyle{Chicago}
\bibliography{KZ_BIB_120118}
\end{document}